\documentstyle[12pt, amsfonts,color]{article}
\newlength{\abstractwidth}
\renewcommand{\thefootnote}{\fnsymbol{footnote}}
\renewcommand{\thanks}[1]{\footnote{#1}} % Use this for footnotes
\newcommand{\starttext}{ \setcounter{footnote}{0}
\renewcommand{\thefootnote}{\arabic{footnote}}}

\newcommand{\be}{\begin{equation}}
\newcommand{\bea}{\begin{eqnarray}}
\newcommand{\eea}{\end{eqnarray}} \newcommand{\ee}{\end{equation}}
 \newcommand{\<}{\langle}
\renewcommand{\>}{\rangle} \def\ba{\begin{eqnarray}}
\def\ea{\end{eqnarray}}

\def\ho{(\frac{\partial}{\partial t} - \Delta)}

\def\ge{\geq}
\def\le{\leq}

\def\qed{$\Box$}

\def\p{\partial}

\def\[{{\bf [}}
\def\]{{\bf ]}}

\def\sigmaone{\sigma}

\def\ric{\textrm{Ric}}

\begin{document}
\starttext \baselineskip=18pt \setcounter{footnote}{0}
\newtheorem{theorem}{Theorem}
\newtheorem{lemma}{Lemma}
\newtheorem{corollary}{Corollary}
\newtheorem{definition}{Definition}
\newtheorem{conjecture}{Conjecture}
\newtheorem{proposition}{Proposition}

\begin{center}
{\Large \bf PARABOLIC DIMENSIONAL REDUCTIONS OF 11D SUPERGRAVITY
\footnote{Work supported in part by the National Science Foundation under grants DMS-12-66033 and DMS-17-10500.}}
\nonumber \\
\medskip
\centerline{
Teng Fei, Bin Guo, and Duong H. Phong}

\begin{abstract}

{\small
Ansatze are constructed under which the solutions of $11D$ supergravity must be  stationary points of a parabolic flow on a Riemannian manifold $M^{10-p}$. This parabolic flow turns out to be the Ricci flow coupled to a scalar field, a $(3-p)$-form, and a $4$-form. This allows the introduction of techniques from parabolic partial differential equations to the search of solutions to $11D$ supergravity. As a first step, Shi-type estimates and criteria for the long-time existence of the flow are established.
}
\end{abstract}

\end{center}

\bigskip

\baselineskip=15pt

\section{Introduction}
\setcounter{equation}{0}

Ever since $11D$ supergravity was constructed by Cremmer, Julia, and Scherk \cite{CJS}, and even more so after the realization that it is a low-energy effective action of M Theory \cite{HW, T, W}, there has been considerable interest in its solutions. Many have been found through a variety of ansatze (see e.g. \cite{duff, DH, PvN, PW, dW1, BB, MS, GP, englert1982} and references therein). Notable supersymmetric examples include compactifications on Einstein $7$-manifolds \cite{freund1980}, on manifolds with special holonomy \cite{PT}, and multi-membrane solutions \cite{duff1991}. A construction of a class of solutions, starting from a Ricci-flat $8$-manifold which is either compact or complete with faster than quadratic volume growth was proposed recently in \cite{FGP}.

\medskip
Mathematically, the low-energy effective actions and compactifications of string theory have led to the discovery of many deep and unexpected phenomena, beginning with the K\"ahler Ricci-flat compactifications of the heterotic string proposed by Candelas, Horowitz, Strominger, and Witten \cite{CHSW}, and the subsequent discovery of mirror symmetry. The generalization by Hull \cite{H1} and Strominger \cite{S} of the proposal of \cite{CHSW} has now been found to have a very rich mathematical structure as well. The first non-K\"ahler solution of the Hull-Strominger system was found by Fu and Yau \cite{FY} using a geometric construction going back to Calabi and Eckmann \cite{CE, GoPr}, and many more solutions have since been found, including an infinite number of topologically distinct types by Fei, Huang, and Picard \cite{FHP1}, generalizing a geometric construction of Calabi \cite{C} and Gray \cite{G} (see e.g. \cite{F} for more references).
The Hull-Strominger system has also motivated the introduction of many new analytic methods, including flows of $(2,2)$-forms \cite{PPZ1, PPZ2, PPZ3, FHP2}.

\medskip
The main goal of this paper is to begin a more systematic, analytic study of the field equations of $11D$ supergravity than has been available so far.
The field equations of $11D$ supergravity are a system of partial differential equations, more specifically Einstein's equation for a metric coupled to a closed $4$-form $F$ (see (\ref{FieldEqs}) below), so we expect that a full analytic understanding will ultimately require the theory of non-linear partial differential equations. Since space-time is Lorentzian, the equations are hyperbolic. Hyperbolic equations are notoriously difficult in general, and many more tools seem available for elliptic and parabolic systems, such as the Calabi-Yau equation \cite{Y}, the Hermitian-Yang-Mills equation \cite{D,UY}, or the Ricci flow \cite{H,P}.
As a practical first step in the long-term program of finding solutions of $11D$ supergravity by partial differential equations methods, we would like then to identify ansatze by which they can be reduced to an elliptic system.

\medskip
More specifically,
the bosonic fields in the theory are a metric $g_{AB}$ and a closed $4$-form on an $11$-dimensional Lorentzian manifold $M^{11}$, and
the field equations are given by
\footnote{A derivation of these equations and our conventions for forms are provided in the appendix. }
\bea
\label{FieldEqs}
&&
d\star F= {1\over 2} F\wedge F
\\
&&
R_{AB}= {1\over 2}F^2_{AB}-{1\over 6}|F|^2g_{AB}
\eea
While the problem of finding dimensional reductions which are parabolic is only non-trivial when $M^{11}$ is Lorentzian, the resulting reduction process works as well for $M^{11}$ Euclidian, and we can treat both cases simultaneously by introducing a parameter $\sigmaone$ which is defined to be $+1$ when $M^{11}$ is Lorentzian and $-1$ when $M^{11}$ is Euclidian.

\smallskip

Our starting point is to
view the field equations (\ref{FieldEqs}) as stationary points of the following dynamical system
\footnote{A priori the set of stationary points may be larger than the set of solutions of the field equations. This is an issue to be examined separately from the considerations of the present paper. A brief discussion is included in Appendix C.}
\bea
\label{11Dflow}
&&
\frac{\partial F}{\partial t}=-\Box_g F - {\sigma\over 2}d\star(F\wedge F)
\\
&& \label{11Dflow1}
\frac{\partial g_{AB}}{\partial t}=-2R_{AB}+F^2_{AB}-{1\over 3}|F|^2g_{AB}
\eea
where $\Box_g=dd^\dagger+d^\dagger d$ is the Hodge-Laplacian.
When the metric $g_{AB}$ is Lorentzian, this flow is not parabolic, and even its short-time existence is not guaranteed in general. However, we shall find ansatze preserved by the flow, under which the most difficult Lorentz components of the metric are static, and the other components evolve by a parabolic flow. Thus assume that the space time $M^{11}$ is a warped product
$M^{11}=M^{1,p}\times M^{10-p}$ with metric $g$ and $4$-form $F$ of the form
\bea
\label{Ansatz1}
g=e^f \tilde g+\hat g,
\qquad
F= dvol_{\tilde g}\wedge \beta+\Psi.
\eea
Here the metric $\tilde g$ is a Lorentzian metric on $M^{1,p}$ if $g_{AB}$ is Lorentzian, $\hat g$ is a Riemannian matric on $M^{10-p}$, $f$ is a scalar function on $M^{10-p}$, $dvol_{\tilde g}$ is the volume form of $\tilde g$ on $M^{1,p}$, and $\beta$ and $\Psi$ are closed $(3-p)$-forms and $4$-forms on $M^{10-p}$ respectively. The metric $\tilde g$ is Riemannian if $g_{AB}$ is Riemannian. The dimension $p$ can take any integer value between $0$ and $10$ (when $p\geq 4$, the form $\beta$ is $0$). Assume that the metric $\tilde g$ is Einstein with scalar curvature $ (1+p)\tilde \lambda$, i.e.,
\bea\nonumber
\ric(\tilde g)=\tilde \lambda \tilde g.
\eea
%Throughout the paper, we will define a number $\sigma$, which is $1$ when the metric $g$ on $M^{11}$ is Lorentzian, and $-1$ when $g$ is Riemannian.

\noindent Then we have the following theorem:

\begin{theorem}
\label{eqn:system}
Consider the following flow of the tuple $(\hat g , f,\beta, \Psi)$ on $M^{10-p}$ with the initial values $\beta_0$ and $\Psi_0$ being closed forms
\bea
\label{eqn:g1 main}
\frac{\partial \hat g_{ij}}{\partial t}&=&-2{\mathrm{Ric}}(\hat g)_{ij}+(p+1)\left((\nabla_{\hat g}^2f)_{ij}+{1\over 2}f_if_j\right)\nonumber\\
&&\quad -\sigmaone e^{(p+1)f}\beta^2_{ij}
+\Psi^2_{ij}-{1\over 3}(|\Psi|^2-\sigmaone e^{-(p+1)f}|\beta|^2)\hat g_{ij}\\
\label{eqn:f1 main}
\frac{\partial f}{\partial t}&=&\Delta_{\hat g} f+{p+1\over 2}
|\nabla_{\hat g} f|_{\hat g}^2-{2\over 3}\sigmaone e^{-(p+1)f}|\beta|_{\hat g}^2
-{1\over 3}|\Psi|_{\hat g}^2-2\tilde\lambda e^{-f}\\
\label{eqn:beta1 main}
\frac{\partial \beta}{\partial t}&=&-\Box_{\hat g} \beta+(-1)^{p+1}{p+1\over 2} d\star_{\hat g} (df\wedge\star_{\hat g}\beta)
- \sigmaone (-1)^p{p+1\over 4}e^{{p+1\over 2}f}df\wedge \star_{\hat g}(\Psi\wedge\Psi)
\nonumber\\
&&\quad- \sigmaone (-1)^p{1\over 2}e^{{p+1\over 2}f}d\star_{\hat g}(\Psi\wedge\Psi)
\\
\label{eqn:psi1 main}
\frac{\partial \Psi}{\partial t} &=&-\Box_{\hat g}\Psi+(-1)^p{p+1\over 2}d\star_{\hat g}(df\wedge\star_{\hat g}\Psi)
 - {p+1\over 2}e^{-{p+1\over 2}f}df\wedge\star_{\hat g}(\beta\wedge\Psi)
\nonumber\\
&&
\quad +e^{-{p+1\over 2}f} d\star_{\hat g}(\beta\wedge\Psi)
\eea
Then the following hold:

{\rm (a)} The forms $\beta$ and $\Psi$ remain closed along the flow, and the pair
$(g, F)$ defined by (\ref{Ansatz1}) satisfies the flow (\ref{11Dflow}) and (\ref{11Dflow1}) on the $11$-dimensional Lorentzian/Riemannian manifold $M^{11}$.

{\rm (b)} Assume that $M^{10-p}$ is compact. Then the above flow (\ref{eqn:g1 main}) - (\ref{eqn:psi1 main}) is weakly parabolic, and admits a smooth solution at least on some interval $[0,T_0)$ for $T_0>0$ depending only on the initial values.

{\rm (c)} If $T<\infty$ is the maximum existence time of the above flow, then
\bea\nonumber
\limsup_{t\to T^-} \sup_{M^{10-p}} ( |Rm| + |f|+|\beta| +  |\Psi|   ) = \infty.
\eea
\end{theorem}

The property (c) implies that the flow is well-behaved, in the sense that, in order for it to terminate, one of only 4 quantities must blow up. Its proof requires estimates of higher order for all fields $g,f,\beta,\Psi$, which are described in section \S 4 below. For the sake of simplicity, we have restricted ourselves here to the case of $M^{10-p}$ compact. But in view of the considerations explained in \cite{FGP}, the case of $M^{10-p}$ non-compact is also of interest, and we shall return to it elsewhere.

\medskip
Supergravity in $11D$ is formulated with a Lorentz signature, but it is likely that a version with Euclidean signature would be of mathematical interest as well, just as Yang-Mills theory with Euclidean signature achieved a prominent position in both geometry and particle physics. For this Euclidean version, we can obtain criteria for long-time existence of the corresponding flow, in analogy with Theorem \ref{eqn:system}:

\begin{theorem}
\label{Euclidean}
Consider the flow (\ref{11Dflow}) and (\ref{11Dflow1}) with $\sigma= -1$ where $g_{AB}$ and $F$ are assumed to be a metric and a closed $4$-form on a compact Riemannian $11$-dimensional manifold $M^{11}$. Then the following hold

{\rm (a)} The $4$-form $F$ remains closed as long as the flow exists.

{\rm (b)} The flow exists at least for a time interval $[0,T_0)$ for some $T_0>0$;

{\rm (c)} If $T<\infty$ is the maximum existence time of the flow (\ref{11Dflow}) and (\ref{11Dflow1}) then
\bea\nonumber
\limsup_{t\to T^-} \sup_{M^{11}} ( |Rm| + |F|  ) = \infty.
\eea
\end{theorem}

The estimates for the higher order derivatives of $Rm $ and $F$ are described in \S 3.

\bigskip
Eleven-dimensional supergravity is a highly constrained physical theory, and arguably geometrically the simplest among the low-energy effective limits of $M$ Theory. As such, it is natural to expect that its solutions may lead to some new canonical structure on eleven-dimensional manifolds and their compactifications.

\medskip
Perhaps not surprisingly, the flows from the dimensional reductions of $11D$ supergravity and of its Euclidean version are all Ricci flows, coupled to tensor fields of different ranks.
Coupled Ricci flows have been considered before in the mathematics literature \cite{Li, L,M, GHP}, in the simplest case of coupling to a single scalar field.
A Perelman pseudo-locality theorem in the case of couplings to a single scalar field was obtained in \cite{GHP}. The coupled flows introduced here will probably require many new techniques, which may also be interesting in their own right from the point of view of the theory of non-linear partial differential equations.

\section{Dimensional Reductions}
\setcounter{equation}{0}

The main goal of this section is to establish part (a) of Theorem \ref{eqn:system}. Parts (b,c) of Theorem \ref{eqn:system} as well as Theorem \ref{Euclidean} require estimates and will be established in subsequent sections.

\medskip
We use indices $A,B,C,\cdots$ for coordinates on $M^{11}$, indices $a,b,c,\cdots$ for coordinates on $M^{1,p}$, and indices $i,j,k,\cdots$ for coordinates on $M^{10-p}$. Then the curvature of the warped product (\ref{Ansatz1}) is given by
\bea
R_{ijkl}=\hat R_{ijkl},
\qquad
R_{iajb}=-e^f\left({1\over 2}(\nabla_{\hat g}^2f)_{ij}+{1\over 4}f_if_j\right)\tilde g_{ab}
\nonumber\\
R_{abcd}=e^f\tilde R_{abcd}+{1\over 4}e^{2f}|\nabla_{\hat g} f|_{\hat g}^2\nonumber
(\tilde g_{ad}\tilde g_{bc}-\tilde g_{ac}\tilde g_{bd}).
\eea
Contracting gives the relation between the Ricci curvatures,
\bea\nonumber
\ric(g)_{ab}=\ric(\tilde g)_{ab}-{1\over 2}e^f(\Delta_{\hat g} f+{p+1\over 2}|\nabla_{\hat g} f|_{\hat g}^2)\tilde g_{ab}
\\
\ric(g)_{ij}=\ric(\hat g)_{ij}-{p+1\over 2}\left((\nabla_{\hat g}^2f)_{ij}+{1\over 2}f_if_j\right).\nonumber
\eea
All the other Ricci curvature components vanish. i.e. $\ric(g)_{ia}=0$, for all $i, a$.

\medskip
Using these formulas, part (a) of Theorem \ref{eqn:system} can now be proved by a direct calculation. However, it may be more instructive to proceed in a way that shows how the dimensionally reduced flows (\ref{eqn:g1 main}) -  (\ref{eqn:psi1 main}) arise.
While the space-time $M^{11}$ is always taken to be the warped product (\ref{Ansatz1}), we shall consider different ansatze for the $4$-form $F$.

\subsection{Possible ansatze}

Ideally we would like our ansatze for $F$ to be as general as possible. Taking the product structure into account, a very general ansatz for $F$ is the following
$$F=\alpha_4+\alpha_3\wedge\beta_1+\alpha_2\wedge\beta_2+\alpha_1\wedge\beta_3+\Psi,$$
where $\alpha_j$ are $j$-forms on $M^{1,p}$, $\beta_j$ are $j$-forms on $M^{10-p}$, and $\Psi$ a $4$-form on $M^{10-p}$. As $F$ is closed, we want all these forms to be closed. Under such an ansatz and if $\tilde{g}$ is Einstein, the $ab$-component of the equation (1.2) becomes
\bea
&&\left(\tilde{\lambda}-{1\over 2}e^f(\Delta_{\hat g} f+{p+1\over 2}|\nabla_{\hat g}f|_{\hat g}^2)+\frac{e^f}{6}|F|_{\hat g}^2\right)\tilde g_{ab}\nonumber\\
&=& \frac{1}{2}\left(e^{-3f}(\alpha_4^2)_{ab} + e^{-2f}|\beta_1|^2(\alpha_3^2)_{ab} + e^{-f}|\beta_2|^2(\alpha_2^2)_{ab} + |\beta_3|^2(\alpha_1^2)_{ab}\right).\nonumber
\eea

\noindent Therefore, a natural assumption for the flow to reduce is that all $\alpha_j$'s satisfy
$$(\alpha_j^2)_{ab}:=\langle \iota_{\p_a}\alpha_j,\iota_{\p_b}\alpha_j\rangle=c_j\cdot\tilde{g}_{ab}$$
for some constant $c_j$. % depending only on $M^{10-p}$ but independent of $M^{1,p}$.

In Riemannian geometry, there are many such differential forms, such as all the calibration forms \cite{HL}, which include powers of the K\"ahler form and the fundamental 3 and 4 forms in $G_2$ geometry. However, such forms are rare in Lorentzian geometry due to the following linear algebra fact

\begin{proposition}\label{form}
Let $(V,\tilde{g})$ be a finite-dimensional Lorentzian vector space. Suppose $\alpha\in\wedge^j V^*$ is a $j$-form such that $\alpha^2_{ab}=c\cdot\tilde{g}_{ab}$ for some constant $c$. Then $\alpha$ must be a constant or a constant multiple of volume form.
\end{proposition}
\noindent{\it Proof. } By choosing an orthonormal basis, we may write
$$\tilde{g}=-(dx^0)^2+(dx^1)^2+\dots+(dx^p)^2,$$
where $\dim V=p+1$. Write
$$\alpha=\frac{1}{j!}\alpha_{i_1\dots i_j}dx^{i_1}\wedge\dots\wedge dx^{i_j}=\sum_{i_1<\dots<i_j}\alpha_{i_1\dots i_j}dx^{i_1}\wedge\dots\wedge dx^{i_j},$$
we only need to show that $j=0$ or $j=p+1$. It is easy to see that
\be\label{alpha0}
\alpha^2_{00}=\sum_{1\leq i_2<\dots<i_j\leq p}\alpha_{0i_2\dots i_j}^2=-c
\ee
and
\be\label{alphak}
\alpha^2_{kk}=\sum_{1\leq i_2<\dots<i_j\leq p}\alpha_{ki_2\dots i_j}^2-\sum_{1\leq i_3<\dots<i_j\leq p}\alpha^2_{0ki_3\dots i_j}=c
\ee
for $k=1,2,\dots,p.$ Summing (\ref{alphak}) over $k$ and add $p$ times of (\ref{alpha0}), we get
$$j\sum_{1\leq i_1<\dots<i_j\leq p}\alpha_{i_1\dots i_j}^2+(p+1-j)\sum_{1\leq i_2<\dots<i_j\leq p}\alpha^2_{0i_1\dots i_j}=0.$$
If $0<j<p+1$, we conclude from above equation that $\alpha=0$, thus the proposition is proved. \hfill\qed

\medskip
Due to Proposition \ref{form}, the most general ansatz we shall consider takes the form
$$F=dvol_{\tilde{g}}\wedge\beta+\Psi,$$
where $\beta$ is a closed $(3-p)$-form and $\Psi$ a closed $4$-form on $M^{10-p}$.

%\smallskip
%The reason for the ansatze below is that we would like space-time to remain a warped product under the dimensionally reduced flow. A priori the $4$-form $F$
%can involve forms on $M^{1,p}$ of different degrees. But the condition that space-time remain a warped product implies that the forms on $M^{1,p}$
%must be either the volume form on $M^{1,p}$ or $0$ (see Appendix C
%for the precise statement and proof).
%%%%%%%%%%%%%%%%%%%%%%
\subsection{The case $F=dvol_{\tilde g}\wedge \beta$}\label{section2.2}

Here for simplicity and as a warm-up, we take $$F = dvol_{\tilde g} \wedge \beta,$$
where $\beta$ is a closed $(3-p)$-form on $M^{10-p}$. We have the following formulae.
$$\star_g F = e^{-\frac{p+1}{2}f} \star_{\tilde g}dvol_{\tilde g} \wedge \star_{\hat g} \beta = - \sigmaone  e^{-\frac{p+1}{2}f} \star_{\hat g} \beta,$$
$$d\star_g F = \sigmaone \frac{p+1}{2} e^{-\frac{p+1}{2}f} df\wedge \star_{\hat g} \beta  - \sigmaone e^{-\frac{p+1}{2}f} d\star_{\hat g} \beta,$$
$$\star_g d \star_g F= - \sigmaone dvol_{\tilde g} \wedge \star_{\hat g} d \star _{\hat g} \beta + \frac{p+1}{2}\sigmaone dvol_{\tilde g} \wedge \star_{\hat g}( df \wedge \star_{\hat g} \beta  ),$$
and
$$d\star_{g} d \star_g F = -\sigmaone (-1)^{p+1} dvol_{\tilde g} \wedge d\star_{\hat g} d \star _{\hat g} \beta + \sigmaone (-1)^{p+1}\frac{p+1}{2} dvol_{\tilde g} \wedge d\star_{\hat g}( df \wedge \star_{\hat g} \beta  ),$$
$$F\wedge F = 0.$$
From $$\frac{\partial F}{\partial t} = dvol_{\tilde g} \wedge \frac{\partial \beta}{\partial t}$$
and the equation (\ref{11Dflow}) we get
$$\frac{\partial \beta}{\partial t} = -  (-1)^{p+1} d\star_{\hat g} d \star_{\hat g} \beta +  (-1)^{p+1} \frac{p+1}{2} d\star_{\hat g} ( df\wedge \star_{\hat g}\beta  ),$$ since $d^\dagger_{\hat g} \beta = (-1)^{p+1} \star_{\hat g} d \star_{\hat g} \beta$, it follows that
\bea\label{eqn:beta form}
\frac{\partial \beta}{\partial t} = - d d^\dagger_{\hat g} \beta + (-1)^{p+1} \frac{p+1}{2} d\star_{\hat g} ( df \wedge \star_{\hat g} \beta  ),
\eea
By direct calculations, we have
$$(F^2)_{ab} = - \sigmaone e^{-p f} |\beta|_{\hat g}^2 \tilde g_{ab}, $$ %\quad{\mathrm{so\,\,}} h = - e^{-pf} |\beta|_{\hat g}^2,$$
and $$(F^2)_{ij} = - \sigmaone e^{- (p+1) f} (\beta^2)_{ij},\quad |F|_g^2 = -\sigmaone e^{-(p+1) f} |\beta|^2_{\hat g}.$$
Thus the equations (\ref{11Dflow}) and (\ref{11Dflow1}) are reduced to ones of $(\hat g, f,\beta)$ on $M^{10-p}$ as follows
\bea\label{eqn:hat g}
\frac{\partial \hat g_{ij}}{\partial t} = - 2 \ric(\hat g)_{ij} + (p+1) ( (\nabla^2_{\hat g} f)_{ij} + \frac 1 2 f_i f_j )- \sigmaone e^{-(p+1) f} (\beta^2)_{ij} + \frac 13 \sigmaone e^{-(p+1) f} |\beta|_{\hat g}^2 \hat g_{ij},
\eea
\bea\label{eqn:f function}
\frac{\partial f}{\partial t} = \Delta_{\hat g} f + \frac{p+1}{2} |\nabla_{\hat g} f |_{\hat g}^2 - \frac 2 3 \sigmaone e^{-(p+1)f} |\beta|_{\hat g}^2 - 2 \tilde \lambda e^{-f}.
\eea
\bea\label{eqn:beta form 1}
\frac{\partial \beta}{\partial t} = - \Box_{\hat g} \beta + (-1)^{p+1} \frac{p+1}{2} d\star_{\hat g} ( df \wedge \star_{\hat g} \beta  ),
\eea
where $\Box_{\hat g} = dd^\dagger_{\hat g} + d^\dagger_{\hat g} d$ is the Hodge Laplacian with respect to $\hat g$.

{\subsection{The case $F=\Psi$}\label{section2.3}}

Next we take
$$F = \Psi,\quad dF = 0,\quad\Psi\in\wedge^4 M^{10-p}.$$
$$\star_g F = \sigmaone e^{\frac{p+1}{2}f} dvol_{\tilde g} \wedge \star_{\hat g} \Psi$$
$$d\star_g F =\sigmaone (-1)^{p+1}\frac{p+1}{2} e^{\frac{p+1}{2}f} dvol_{\tilde g}\wedge  df\wedge \star_{\hat g}\Psi + \sigmaone(-1)^{p+1} e^{\frac{p+1}{2}f} dvol_{\tilde g} \wedge d\star_{\hat g}\Psi.$$
$$\star_g d\star_g F =\sigmaone (-1)^p \frac{p+1}{2}  \star_{\hat g}( df\wedge \star_{\hat g}\Psi  ) +\sigmaone (-1)^p  \star_{\hat g} d\star_{\hat g} \Psi.  $$
$$d\star_g d\star_g F = \sigmaone(-1)^p \frac{p+1}{2}  d \star_{\hat g}( df\wedge \star_{\hat g}\Psi  ) +\sigmaone (-1)^p  d\star_{\hat g} d\star_{\hat g} \Psi.$$
If $10-p\le 7$, then $F\wedge F = 0$; if $10-p\ge 8$, $F\wedge F \neq 0$, then
$$\star_g(F\wedge F) =  e^{\frac{p+1}{2}f} dvol_{\tilde g} \wedge \star_{\hat g}(\Psi\wedge \Psi),$$
$$d\star_g (F\wedge F) =   (-1)^{p+1}\frac{p+1}{2} e^{\frac{p+1}{2}f} dvol_{\tilde g}\wedge df \wedge \star_{\hat g}(\Psi\wedge \Psi)  +  (-1)^{p+1} e^{\frac{p+1}{2} f}dvol_{\tilde g} \wedge d\star_{\hat g}(\Psi\wedge \Psi)$$
Therefore, in case $10-p\ge 8$, the condition that $F$ remains a closed $4$-form $\Psi$ on $M^{10-p}$  in general may not be preserved along the flow. So here we assume $10-p\le 7$ and for dimension reason $F\wedge F = 0$, and the equation (\ref{11Dflow}) for $F$ becomes
$$\frac{\partial \Psi}{\partial t} = (-1)^p d\star_{\hat g} d\star_{\hat g}\Psi + (-1)^p \frac{p+1}{2} d\star_{\hat g}(df\wedge \star_{\hat g} \Psi),   $$
since $d^\dagger_{\hat g}\Psi = - (-1)^p \star_{\hat g}d\star_{\hat g}\Psi$, the equation above is
$$\frac{\partial \Psi}{\partial t} = -d d^\dagger_{\hat g}\Psi + (-1)^p \frac{p+1}{2} d\star_{\hat g}(df\wedge \star_{\hat g} \Psi),   $$
Straightforward calculations show that
$$(F^2)_{ab} = 0,\quad (F^2)_{ij} = (\Psi^2)_{ij},\quad |F|^2_{g} = |\Psi|^2_{\hat g}$$
Hence the flow equations (\ref{11Dflow}) and (\ref{11Dflow1}) are reduced to the following equations on the Riemannian part $M^{10-p}$ with $10-p\le 7$
\bea\label{eqn:g 3}
\frac{\partial \hat g_{ij}}{\partial t} = - 2 \ric(\hat g)_{ij} + (p+1) ( (\nabla^2_{\hat g} f)_{ij} + \frac 1 2 f_i f_j  ) + (\Psi^2)_{ij} - \frac 1 3 |\Psi|^2 \hat g_{ij}.
\eea
\bea\label{eqn:f 3}
\frac{\partial f}{\partial t} =\Delta_{\hat g} f + \frac{p+1}{2} |\nabla_{\hat g} f|_{\hat g}^2 - \frac{1}{3}|\Psi|^2_{\hat g} - 2 \tilde\lambda e^{-f},
\eea
\bea\label{eqn:F 3}
\frac{\partial \Psi}{\partial t} = -\Box_{\hat g}\Psi + (-1)^p \frac{p+1}{2} d\star_{\hat g}(df\wedge \star_{\hat g} \Psi).
\eea

%%%%%%%%%%%%%%%%%%%%%%%%%%%%%%%%%%

\subsection{The case $F=dvol_{\tilde g}\wedge \beta+\Psi$}\label{section2.4}

In this case, we assume $\beta = 0$ if ${\rm dim}\,M^{1,p} = 1+p\ge 5$. For dimensional reasons we also have $\Psi\wedge\Psi = 0$ in the equations below if ${\rm dim}\,M^{10-p} = 10-p\le 7$.

\medskip

\noindent Let us now take $$F =  dvol_{\tilde g} \wedge \beta + \Psi,$$ where $\beta$ and $\Psi$ are chosen as in the previous two cases.
Similar to the previous equations we have
\bea\nonumber
d\star_g d\star_g F & = &  - \sigmaone (-1)^{p+1} dvol_{\tilde g} \wedge d\star_{\hat g} d \star _{\hat g} \beta + \sigmaone(-1)^{p+1}\frac{p+1}{2} dvol_{\tilde g} \wedge d\star_{\hat g}( df \wedge \star_{\hat g} \beta  )\\
&& \nonumber+\sigmaone (-1)^p \frac{p+1}{2}  d \star_{\hat g}( df\wedge \star_{\hat g}\Psi  ) +\sigmaone (-1)^p  d\star_{\hat g} d\star_{\hat g} \Psi.
\eea
From
$$F\wedge F = 2 dvol_{\tilde g}\wedge \beta\wedge \Psi + \Psi\wedge \Psi$$
we get
\bea\nonumber
d\star_g (F\wedge F)& = &  \sigmaone (p+1)e^{-\frac{p+1}{2}f} df\wedge \star_{\hat g}\Psi - 2\sigmaone e^{-\frac{p+1}{2}f}d\star_{\hat g} \Psi\\
&& \nonumber + (-1)^p\frac{p+1}{2}e^{\frac{p+1}{2}f} dvol_{\tilde g}\wedge df \wedge \star_{\hat g}(\Psi\wedge \Psi) + (-1)^p e^{\frac{p+1}{2} f}dvol_{\tilde g} \wedge d\star_{\hat g}(\Psi\wedge \Psi)
\eea
From the equation (\ref{11Dflow}) on $F$, we can derive that
\bea\label{eqn:beta 1}
\frac{\partial \beta}{\partial t} & = & - (-1)^{p+1} d\star_{\hat g} d \star_{\hat g} \beta + (-1)^{p+1} \frac{p+1}{2} d\star_{\hat g}(df \wedge \star_{\hat g} \beta)\\
&& \nonumber- \sigmaone  \frac{(-1)^p (p+1)}{4} e^{\frac{p+1}{2}f} df\wedge \star_{\hat g}(\Psi\wedge \Psi)- \sigmaone  \frac{(-1)^p}{2} e^{\frac{p+1}{2}f} d\star_{\hat g}(\Psi\wedge \Psi),\nonumber
\eea
and
\bea\label{eqn:psi 1}
\frac{\partial \Psi}{\partial t} & = & (-1)^p d\star_{\hat g} d \star_{\hat g} \Psi + (-1)^p \frac{p+1}{2} d\star_{\hat g} (df \wedge \star_{\hat g} \Psi) \\
&& \nonumber -  \frac{p+1}{2} e^{-\frac{p+1}{2}f} df\wedge \star_{\hat g} \Psi + e^{-\frac{p+1}{2}f} d\star_{\hat g}\Psi.
\eea
We also have
$$(F^2)_{ab} = -\sigmaone e^{-pf} |\beta|_{\hat g}^2 \tilde g_{ab},\quad (F^2)_{ij} = - \sigmaone e^{-(p+1)f} (\beta^2)_{ij} + (\Psi^2)_{ij}$$
$$|F|^2 = -\sigmaone e^{-(p+1)f} |\beta|_{\hat g}^2 + |\Psi|_{\hat g}^2.$$
Then the equations (\ref{11Dflow}) and (\ref{11Dflow1}) are reduced to $M^{10-p}$ as follows:
\bea\label{eqn:hat g 4}
\frac{\partial \hat g_{ij}}{\partial t} &=& - 2 \ric(\hat g)_{ij} + (p+1)((\nabla^2 _{\hat g}f)_{ij} + \frac 12 f_i f_j  ) - \sigmaone e^{-(p+1)f} (\beta^2)_{ij} \\
&&\nonumber + (\Psi^2)_{ij} -\frac 1 3( -\sigmaone e^{-(p+1)f} |\beta|_{\hat g}^2 + |\Psi|_{\hat g}^2  ) \hat g_{ij}
\eea
\bea\label{eqn:f 4}
\frac{\partial f}{\partial t} & = & \Delta_{\hat g} f + \frac{p+1}{2} |\nabla_{\hat g} f|_{\hat g}^2 - \frac{2}{3}\sigmaone e^{-(p+1)f}(\beta^2)_{ij} -\frac 1 3 |\Psi|^2_{\hat g} - 2 \tilde \lambda e^{-f}.
\eea
and combining with equations (\ref{eqn:beta 1}) and (\ref{eqn:psi 1}) and the assumption $d\beta = 0$ and $d\Psi = 0$ we derive the evolution equations for $\beta$ and $\Psi$, which are the main equations (\ref{eqn:beta1 main}) and (\ref{eqn:psi1 main}) we are going to study:
\bea\label{eqn:beta new}
\frac{\partial \beta}{\partial t} &= & - \Box_{\hat g} \beta + (-1)^{p+1} \frac{p+1}{2} d\star_{\hat g}(df \wedge \star_{\hat g} \beta)\\
&& \nonumber- \sigmaone (-1)^p \frac{p+1}{4} e^{\frac{p+1}{2}f} df\wedge \star_{\hat g}(\Psi\wedge \Psi)- \sigmaone  \frac{(-1)^p}{2} e^{\frac{p+1}{2}f} d\star_{\hat g}(\Psi\wedge \Psi),\nonumber
\eea
and
\bea\label{eqn:psi new}
\frac{\partial \Psi}{\partial t} &= & - \Box_{\hat g} \Psi + (-1)^p \frac{p+1}{2} d\star_{\hat g} (df \wedge \star_{\hat g} \Psi) \\
&& \nonumber - \frac{p+1}{2} e^{-\frac{p+1}{2}f} df\wedge \star_{\hat g} \Psi + e^{-\frac{p+1}{2}f} d\star_{\hat g}\Psi.
\eea

\bigskip
This case in section \ref{section2.4} is the most general of all three cases in sections \ref{section2.2}, \ref{section2.3} and \ref{section2.4}, and it is the one considered in Theorem \ref{eqn:system}.

%%%%%%%%%%%%%%%%%%%%%%%%%%

\section{Estimates for the $11D$ Euclidean Flow}
\setcounter{equation}{0}
In this section, we take $\sigma = -1$ in (\ref{11Dflow}).
We prove first Theorem \ref{Euclidean}, as the flow (\ref{11Dflow}) and (\ref{11Dflow1}) is technically simpler than the dimensionally reduced flow (\ref{eqn:g1 main}) - (\ref{eqn:psi1 main}) considered in Theorem \ref{eqn:system}. Since in this section we consider general compact Riemannian manifolds $M^{11}$ and not any particular warped product, we shall lighten the notation and use $i,j,k,\cdots$ to index coordinates on $M^{11}$   instead of $A,B,C,\cdots$ as in the other sections.

\medskip
Part (a) of Theorem \ref{Euclidean} follows easily from the
fact that $d\Box_g = \Box_g d$ and hence $\partial_t ( d F ) = - \Box_g (dF)$.
In particular, if $dF$ is $0$ at time $t=0$, then $dF(t)$ is $0$ for all time $t>0$ as long as  the solution exists.

\subsection{The short-time existence of the flow}

Next we establish part (b).
We adapt DeTurck's trick to make the equations (\ref{11Dflow}) and (\ref{11Dflow1}) strictly parabolic after re-parameterizations. We consider the following flow for a metric $\hat g_{ij}$ (not to be confused with the metric $\hat g_{ij}$ in warped products considered earlier) and a closed $4$-form $\hat F$ on $M^{11}$,
\bea
\label{eqn:main2}
\frac{\partial \hat F}{\partial t} &= &- \Box_{\hat  g} \hat F + \frac 1 2 d \big(\star (\hat F\wedge \hat F)\big) + L_V \hat F\nonumber \\
\frac{\partial \hat g_{ij}}{\partial t} & =& - 2 \ric(\hat g)_{ij}+ \hat \nabla_i V_j + \hat \nabla_j V_i +   + (\hat F^2)_{ij} - \frac 1 3 |\hat F|^2 \hat g_{ij},
\eea
where the vector field $V = V_i \frac{\partial}{\partial x^i}$ is given by
$$V_i  = \hat g_{i j} \hat g^{kl} (\hat \Gamma^j_{kl} - \Gamma(g_0)^j_{kl} )$$
$\hat \Gamma^j_{kl} = \frac{1}{2} \hat g^{jp}(\partial_{k} \hat g_{pl} + \partial_{l} \hat g_{p k} - \partial_{p} \hat g_{kl }   )$ is the Christoffel symbol of $\hat g$, and similar definition for $\Gamma(g_0)$ and $g_0$ is the initial metric. $L_V \hat F$ denotes the Lie derivative of the $4$-form $\hat F$ in the direction of $V$.

\begin{lemma}\label{lemma2}
The operators $\hat g_{ij}\mapsto -2 \ric(\hat g)_{ij} + \hat \nabla_i V_j + \hat \nabla_j V_i$ and $\hat F \mapsto - \Box_{\hat  g} \hat F$ are both strictly elliptic.
\end{lemma}
\noindent {\it Proof.} This lemma is well-known as the DeTurck's trick in the study of Ricci flow. For the convenience of  readers, we provide a proof below.

The ellipticity of the operator on $\hat g_{ij}$ follows from a straightforward calculation. At any fixed point $p\in M$, we choose normal coordinates for $g_0$ so that $(g_0)(p)_{ij} = \delta_{ij}$ and $dg_0(p) = 0$. Then at $p$ we have
\bea
&& \hat \nabla_i V_j + \hat \nabla_j V_i\nonumber \\
&= & \partial_i \hat g_{jk} \hat g^{pq} \hat \Gamma^k_{pq} + \partial_j \hat g_{ik} \hat g^{pq}\hat \Gamma^k_{pq} + \hat g_{jk} \partial_i \hat g^{pq} \hat\Gamma_{pq}^k + \hat g_{ik}\partial_j \hat g^{pq} \hat \Gamma^k_{pq}\nonumber\\
& &+ \hat g_{jk} \hat g^{pq} \partial_i \hat \Gamma^k_{pq} + \hat g_{ik} \hat g^{pq} \partial_j \hat\Gamma^k_{pq} - \hat g_{jk}\hat g^{pq} \partial_i \Gamma(g_0)_{pq}^k  -  \hat g_{ik}\hat g^{pq}\partial_j \Gamma(g_0)^k_{pq} - 2 \hat g_{pk} \hat g^{ts} \hat \Gamma^k_{ts} \hat \Gamma^p_{ij},\nonumber
\eea
and
$$\ric(\hat g)_{ij} = \partial_k \hat\Gamma^k_{ij} - \partial_j \hat\Gamma^k_{ki} + \hat\Gamma^p_{ji} \hat\Gamma^k_{pk} - \hat\Gamma^p_{ki} \hat\Gamma^k_{pj},$$
so the leading (i.e. second) order terms in $-2 \ric(\hat g)_{ij} + \hat \nabla_i V_j + \hat\nabla_j V_i$ are given by
\bea
& &-2 \partial_k \hat\Gamma^k_{ij} + 2 \partial_j \hat\Gamma^{k}_{ki} + \hat g_{jk}\hat g^{pq} \partial_i \hat\Gamma^k_{pq} + \hat g_{ik}\hat g^{pq} \partial_j \hat\Gamma^k_{pq}\nonumber \\
&= & \hat g^{kl} \frac{\partial^2 \hat g_{ij}} {\partial x^k\partial x^l } - 2 \partial_k \hat g^{kl} \hat\Gamma^k_{ij} \hat g_{kl} + 2 \partial_j \hat g^{kl} \hat\Gamma ^p_{ki} \hat g_{pl} + \hat g_{jk} \hat g^{pq} \partial_i \hat g^{kl} \hat\Gamma^m _{pq} \hat g_{ml} + \hat g_{ik} \hat g^{pq} \partial_j \hat g^{kl} \hat\Gamma^m_{pq} \hat g_{m l},\nonumber
\eea
establishing our claim.

As for the operator on $\hat F$, it is the negative Hodge-Laplacian, and its strict ellipticity is well-known. More explicitly, by the Lichnerowicz-Weitzenb\"ock formula, we have
\bea
\label{DeltaF}
 - \Box_{\hat  g} \hat F=\hat g^{pq}\nabla_p \nabla_q  \hat F + \hat  F * \widehat{Rm}
 \eea
where $\hat F* \widehat{Rm}$ denotes the terms which are pointwise linear combination of the product of $\hat F$ and $\widehat{Rm}$. The strict ellipticity of $ - \Box_{\hat g}$ follows. \hfill \qed

\medskip

Lemmas \ref{lemma2} implies that the system (\ref{eqn:main2}) is a system of strictly parabolic differential equations, so the short time existence of $(\hat g(t),\hat F(t))$ for some time interval $[0,T_0)$ with $T_0>0$ follows from standard theory of parabolic partial differential equations.

\medskip
Let $\varphi_t: M\to M$ be  the  one-parameter subgroup of diffeomorphisms generated by $- V$, and define $$g(t) = \varphi_t^* \hat g(t) ,\quad F(t) = \varphi_t^* \hat F(t).$$
It is straightforward to check that $(g(t),F(t))$ solves the equations (\ref{11Dflow}) and (\ref{11Dflow1}).
This proves part (b) of Theorem \ref{Euclidean}.

\subsection{The long-time existence of the flow}

We now establish part (c) of Theorem \ref{Euclidean}, which is that
the flow exists as long as the $|Rm|$ and $|F|^2$ remain bounded. For this, we will prove that boundedness of $|Rm|$ and $|F|^2$ implies that  all the higher order derivatives of $(g, F)$ remain bounded along the flow  (\ref{11Dflow}) and (\ref{11Dflow1}). This type of estimates have been studied extensively for other geometric flows like Ricci flow or Anomaly flow(for example \cite{H, Sh, Li, PPZ2}).

\medskip
We will use the constant $C$ to denote a uniform constant depending only on the dimension $n=11$, which may vary from line to line.
We list the evolution equations of some key quantities. They can be obtained by a direct calculation, so we omit the proof.

\begin{lemma}
Along the flow  (\ref{11Dflow}) and (\ref{11Dflow1}),  $|F|^2$ satisfies
\bea
\frac{\partial}{\partial t} |F|^2 & =&  \Delta_g |F|^2 - 2 |\nabla F|^2 + Rm*F*F + \frac 1{2\cdot 4!} F^{ijkl} \big(d\star(F\wedge F)\big)_{ijkl}  \nonumber \\
&  &+ \Big(2 R_{ij} - (F^2)_{ij} + \frac{1}{3} |F|^2 g_{ij}  \Big) (F^2)_{ij}\nonumber \\
&\le &\Delta |F|^2  - 2 |\nabla F|^2 + C |Rm| |F|^2 + c_0|\nabla F| |F|^2 - \frac{4}{33} |F|^4,\nonumber
\eea
for some uniform constants $C>0$ and $c_0\ge 0$.
\end{lemma}
In the last inequality we apply the Cauchy-Schwarz  inequality $F^2_{ij} F^2_{ij}\ge \frac{(g^{ij} F^2_{ij})^2}{11} = \frac{16}{11} (|F|^2)^2$. We  remark that the term involving $c_0$ comes from $F^{ijkl}(d\star (F\wedge F))_{ijkl}$, and if this number $c_0$ is small, say $c_0^2 \le \frac{32}{33}$, then the long time existence of the flow can be weakened to the condition that $|Rm|$ remains bounded along the flow, since the bound on $|F|^2$ can be derive from the equation above by maximum principle. In particular, if for some choice of $F$, $F\wedge F = 0$ along the flow (e.g. in the dimension reduction of the flow), this observation may be applied.

\begin{lemma}\label{lemma3 3}
The evolution equations for the curvatures are given by
\bea
\frac{\partial }{\partial t} R_{ijkl} & =& \Delta R_{ijkl} + 2 R_{mlnk} R_{mjni} - 2 R_{mknl} R_{mjni} + 2 R_{mkni} R_{mlnj} - 2 R_{mlni} R_{mknj} \nonumber \\
&& - R_{mi} R_{mjkl} - R_{mj} R_{imkl} - R_{mk} R_{ijml} - R_{ml} R_{ijkm}\nonumber\nonumber \\
&& + \frac{1}{2} (F^2)_{ml} R_{mkji} + \frac{3}{2} (F^2)_{km} R_{mlji} - \frac{1}{3} |F|^2 R_{lkji}\nonumber \\
&& + \frac 12 \big((F^2)_{kj,li} - (F^2)_{jl,ki}  -  (F^2)_{kl,ij} - (F^2)_{ki,lj} \big)\nonumber \\
& &+ \frac{1}{6}\big( - g_{kj}(|F|^2)_{,li}  + g_{jl}(|F|^2)_{,ki}  +  g_{ij}(|F|^2)_{,kl} + g_{ki} (|F|^2)_{,lj}   \big)\nonumber \\
&=& \Delta R_{ijkl} + Rm* Rm + F * F * Rm + \nabla^2 F * F + \nabla F * \nabla F,\nonumber
\eea where for a tensor $\alpha$ for simplicity we denote $\alpha_{,i}$ to be the covariant derivative $\nabla_i \alpha$.
%
%\bea
%\frac{\partial }{\partial t} R_{ik} &=&  \Delta R_{ik} - 2 R_{mink}R_{mn} - 2 R_{mi} R_{mk} - \frac{1}{12} (F^2)_{ml} %R_{mkli} - \frac{1}{4} (F^2)_{km} R_{mi} \nonumber \\
%&& + \frac 1{12} \big( (F^2)_{kl,li} - (F^2)_{ll,ki} - (F^2)_{kl,il} - (F^2)_{ki,ll}   \big)\nonumber \\
%&& + \frac{1}{144}\big( -|F|^2_{,li} g_{kl} + 11 |F|^2_{,ik} + |F|^2_{,kl} g_{il} +|F|^2_{,ll} g_{ki}    \big).\nonumber
%\eea
\end{lemma}

%These formulas follow from straightforward calculations, so we omit the proof.
From now on,  for a given  $T>0$, we
assume that there is a finite constant $K>0$ such that
$$\sup_{M\times[0,T)}( |Rm|+ |F|^2) \le K.$$

\begin{lemma}\label{prop:grad} There exists a constant $C_{K,T}$ depending only on $K, T$ such that
$$\sup_M |\nabla F|^2 \le \frac{C_{K,T}}{t},\quad \forall\; t\in (0, T).$$
\end{lemma}

\noindent To prove Lemma \ref{prop:grad}, we need the evolution equation for the quantity $|\nabla F|^2 = \frac{1}{4!}F_{ijkl,p} F^{ijkl,p}$.

\begin{lemma}\label{lemma5 3}
The quantity $|\nabla F|^2$ satisfies
\bea
\frac{\partial}{\partial t} |\nabla F|^2& \le& \Delta |\nabla F|^2 - 2 |\nabla ^2 F|^2 + C\Big( |Rm| |\nabla F|^2 +  |\nabla Rm| |F| |\nabla F|\nonumber \\
&&  +  |\nabla F|^3 +  |F| |\nabla F| |\nabla^2 F| +  |F|^2  |\nabla F|^2\Big),\nonumber
\eea
for some constant $C=C(K)>0$.

\end{lemma}
\noindent{\it Proof:}
By definition
$$F_{ijkl,p} = \partial_p F_{ijkl} - F_{qjkl} \Gamma^q_{pi} - F_{iqkl}\Gamma^q_{pj} - F_{ijql}\Gamma^q_{pk} - F_{ijkq}\Gamma^q_{pl},$$
where $\Gamma^k_{ij} = \frac 1 2 g^{kl}( \partial_i g_{jl} + \partial_j g_{il} - \partial_l g_{ij}  ) $ is the Christoffel symbol of the metric $g$ which satisfies the evolution equation
\bea
\frac{\partial}{\partial t} \Gamma_{ij}^k &= & \frac 1 2 g^{kl} \Big( \nabla_i \big(-2 R_{jl} +  (F^2)_{jl} - \frac{1}{3} |F|^2 g_{jl} \big) + \nabla_j \big( -2 R_{il} +  (F^2)_{il} - \frac 1{3} |F|^2 g_{il}   \big)\nonumber \\
& & - \nabla_l \big( -2 R_{ij} +  (F^2)_{ij} - \frac{1}{3} |F|^2 g_{ij}    \big)   \Big)\nonumber \\
&=& \nabla \ric + F* \nabla F.\nonumber
\eea
Hence it follows that
\bea
\frac{\partial }{\partial t} F_{ijkl,p}&= & \nabla_p \frac{\partial F_{ijkl}}{\partial t} - F * \frac{\partial \Gamma}{\partial t}\nonumber \\
&=& \nabla_p \big( \Delta F_{ijkl} + Rm* F+ F*F   \big) + F* \nabla Ric + F* F*\nabla F\nonumber \\
&=& \nabla_p \Delta F_{ijkl} + \nabla Rm * F + Rm * \nabla F + F * \nabla F + F*F*\nabla F.\nonumber
\eea
Differentiating the formula (\ref{DeltaF}) and commuting derivatives, we obtain
\bea
\nabla_p(-\Box_g F_{ijkl})=  \Delta F_{ijkl,p} + \nabla Rm *F + Rm * \nabla F \nonumber
\eea
thus
$$\frac{\partial }{\partial t}\nabla F = \Delta \nabla F + \nabla Rm * F + Rm * \nabla F + F*\nabla F + F*F*\nabla F.$$
The lemma follows from this equation and the definition $|\nabla F|^2 =\frac{1}{4!} F_{ijkl,p} F^{ijkl,p}$.

\bigskip

\noindent{\it Proof of Lemma \ref{prop:grad}: }
We will use a constant $C_K$ to denote a constant depending only on $K$ and $n=11$, which may vary from line to line.
Set
$$G_1: = t\Big( (|F|^2  + A  )|\nabla F|^2 + |Rm|^2  \Big) + A_1 |F|^2,$$
for some constants $A>0, \, A_1>0$ to be determined.
From the equations in Lemma \ref{lemma3 3} and Lemma \ref{lemma5 3}, we have
\bea
&&(\frac{\partial}{\partial t} - \Delta) \big( (|F|^2 + A) |\nabla F|^2 + |Rm|^2 \big)\nonumber \\
&\le & - |\nabla F|^4 + C_K |\nabla F|^2 - 2(A + |F|^2)|\nabla^2 F|^2 + C_K (A+ |F|^2) |\nabla F|^2 \nonumber \\
&& + C_K (A+|F|^2 ) |\nabla F ||\nabla Rm| + C_K (A+|F|^2 ) |\nabla F|^3 - 2 \langle\nabla |F|^2,\nabla |\nabla F|^2   \rangle\nonumber \\
& &- 2 |\nabla Rm|^2 + C_K |\nabla^2 F| + C_K.\nonumber
 \eea
We observe that by Kato's inequality
$$
|2 \langle\nabla |F|^2,\nabla |\nabla F|^2   \rangle|\le  C_K |\nabla F|^2 |\nabla^2 F|\le \frac 1 {10} |\nabla F|^4 + C_K |\nabla^2 F|^2
$$ and
$$
C_K(A+|F|^2) |\nabla^2 F| |\nabla F| \le C_K^2(A+|F|^2)^2 |\nabla F|^2 + \frac 1 {10} |\nabla^2 F|^2,
$$
$$
C_K(A+|F|^2) |\nabla Rm| |\nabla F| \le C_K^2 (A+|F|^2)^2 |\nabla F|^2 + \frac 1 {10} |\nabla Rm|^2,
$$
$$
C_K|\nabla^2 F|\le \frac 1 {10} |\nabla^2 F|^2 + C_K,
$$and by Young's inequality
$$ \quad C_K (A+|F|^2) |\nabla F|^3\le \frac 1 {10}|\nabla F|^4 + C_K (A+|F|^2)^4.$$
Combining the above inequalities and choosing $A$ sufficiently large (but depending only on $K$) so that the terms involving $|\nabla^2 F|^2$ and $|F|^4$ can be absorbed, we get
\begin{equation}\label{eqn:guo1}
(\frac{\partial}{\partial t} - \Delta) \big( (|F|^2 + A) |\nabla F|^2 + |Rm|^2 \big)
\le - A |\nabla^2 F|^2 - \frac 1 2 |\nabla F|^4 - |\nabla Rm|^2 + C_K.
 \end{equation}
 Therefore if $A_1$ is chosen to be large enough, then
\bea
(\frac{\partial}{\partial t} - \Delta ) G_1 &= & t (\frac{\partial}{\partial t} - \Delta) \big( (|F|^2 + A) |\nabla F|^2 + |Rm|^2 \big) +  \big( (|F|^2 + A) |\nabla F|^2 + |Rm|^2 \big) \nonumber \\
&& + A_1 (\frac{\partial}{\partial t} - \Delta ) |F|^2\nonumber \\
&\le  & t C_K +   (|F|^2 + A) |\nabla F|^2 + |Rm|^2  - A_1 |\nabla F|^2 + A_1 C_K\nonumber \\
&\le & C_{K,T}\nonumber
\eea
The estimate that $G_1\le C_{K,T}$ follows from the maximum principle.
Lemma \ref{prop:grad} is proved from the definition of $G_1$. \hfill \qed

\medskip

Next we have the following higher order derivative estimates.

\begin{lemma}\label{prop:2nd}
There exists a constant $C_{K,T}>0$ such that
$$\sup_M ( |\nabla Rm| + |\nabla ^2 F|  )\le \frac{C_{K,T}}{t},\quad \forall t\in (0, T).$$
\end{lemma}

To prove Lemma \ref{prop:2nd}, we need the following equations which follow from standard calculations so we omit the proof.

\begin{lemma}\label{lemma 7}
We have the following equations:
\bea
\frac{\partial}{\partial t} \nabla^2 F &=&  \Delta \nabla^2 F + \nabla^2 F * Rm + \nabla F * \nabla Rm + F* \nabla^2 Rm + \nabla^2 F * \nabla F\nonumber \\
&& + \nabla^3 F * F + F* \nabla F * \nabla F + F * F * \nabla^2F.\nonumber
\eea
\bea
\frac{\partial}{\partial t}\nabla Rm &=& \Delta \nabla Rm + \nabla Rm * Rm + \nabla Rm * F * F + Rm * F * \nabla F\nonumber \\
& &+ \nabla^3 F * F + \nabla F * \nabla^2 F.\nonumber
\eea
\end{lemma}
%These equations

\bigskip

\noindent{\it Proof of Lemma \ref{prop:2nd}: }
From Lemma \ref{lemma 7} we have the following inequalities:
\bea
\frac{\partial}{\partial t} |\nabla^2 F|^2 & \le & \Delta |\nabla^2 F|^2 - 2 |\nabla^3 F|^2 + C_{K,T}\Big( |\nabla^2 F|^2  +  |\nabla F | |\nabla^2 F| |\nabla Rm| +  |\nabla^2 Rm||\nabla^2 F| \nonumber \\
&&\quad + |\nabla F| |\nabla^2 F|^2 +  |\nabla^2 F| |\nabla^3 F| +  |\nabla F|^2 |\nabla^2 F|\Big),\nonumber
\eea
and
\bea
\frac{\partial}{\partial t}|\nabla Rm|^2& \le & \Delta |\nabla Rm|^2 - 2 |\nabla^2 Rm|^2 + C_K |\nabla Rm|^2 + C_K |\nabla Rm| |\nabla F| \nonumber \\
&& + C_K |\nabla Rm| |\nabla^3 F| + C |\nabla F| |\nabla ^2 F | |\nabla Rm|.\nonumber
\eea
By Cauchy-Schwarz inequality, we have
\bea
& &(\frac{\partial}{\partial t} - \Delta ) \big( |\nabla^2 F| ^2 + |\nabla Rm|^2 \big)\nonumber \\
&\le & - |\nabla^3 F|^2 + C_K\big(1+ |\nabla F| + |\nabla F|^2\big)|\nabla^2 F|^2 + C_K |\nabla F|^2 \nonumber \\
& &-|\nabla ^2 Rm|^2 + C_K |\nabla Rm|^2\nonumber
\eea
We define a quantity
$$G_2 : = t^2 ( |\nabla^2 F|^2 + |\nabla Rm|^2  ) + A_1 t \big((A+|F|^2)|\nabla F|^2 + |Rm|^2\big) + A_2 |F|^2,$$
for the constants $A>0,\, A_1>0, A_2>0$ to be determined. We calculate using the equation (\ref{eqn:guo1}) that
\bea
&& (\frac{\partial }{\partial t} - \Delta )G_2\nonumber \\
&\le & 2 t (|\nabla^2 F|^2 + |\nabla Rm|^2  ) - t^2 |\nabla^3 F|^2 + t^2C_K\big(1+ |\nabla F| + |\nabla F|^2\big)|\nabla^2 F|^2 +t^2 C_K |\nabla F|^2 \nonumber \\
& &-t^2|\nabla ^2 Rm|^2 + C_K t^2|\nabla Rm|^2   + A_1 \big((A+|F|^2) |\nabla F|^2  + |Rm|^2 \big) + C_K A_1 t\nonumber \\
& &-A A_1 t |\nabla^2 F|^2 - \frac 1 2 A_1 t |\nabla F|^4 - A_1 t |\nabla Rm|^2 - A_2 |\nabla F|^2 + C_K\nonumber \\
&\le &  C_{K,T},\nonumber
\eea
here the last inequality is obtained by first choosing $A_1$ large which depends only on $K$ and $T$, then picking $A_2$ big enough. Applying maximum principle to $G_2$, we get $\sup_{M\times[0,T)} G_2\le C_{K,T}$, and this proves Lemma \ref{prop:2nd}.\hfill \qed

\medskip

In general the following higher order estimates can be proved similarly by induction.
\begin{theorem}
\label{theorem3}
There exists a constant $C_{K,T,m}$ for any $m\in\mathbf Z_+$ such that
$$\sup_M ( |\nabla^{m-1} Rm| + |\nabla^m F|  )\le \frac{C_{K,T,m}}{t^{m/2}},\quad \forall~ t\in (0,T).$$
\end{theorem}
\noindent{\it Proof: } In the cases $m=1$ or $m=2$, this inequality is proved in Lemma \ref{prop:grad} and Lemma \ref{prop:2nd}.

\noindent For general $m$, we write
\bea   G_m & :=&  t^m ( |\nabla^m F|^2 + |\nabla^{m-1} Rm|^2 ) \nonumber\\
&&+ \sum_{i=2}^{m-1} A_i t^i (|\nabla^i F|^2  + |\nabla^{i-1} Rm|^2 ) + A_1 t ( (A_0+|F|^2) |\nabla F|^2 +|Rm|^2   ) + B |F|^2,\nonumber
\eea
for suitable choices of constants $A_i>0$ and $B>0$, by an induction argument similar to (but simpler than) the proof of Theorem \ref{prop:4} in next section, we can prove that
$$(\frac{\partial}{\partial t} - \Delta ) G_m\le C_{K,T,m},$$ from which we get $\sup_{M\times [0,T)} G_m\le C_{K,T,m}$, which implies the desired inequality. \hfill\qed

\bigskip
We can now prove part (c) of Theorem \ref{Euclidean}.
We argue by contradiction. If $$\limsup_{t\to T^-} \sup_{M^{11}} ( |Rm| + |F||) =:K <\infty.$$
Then Theorem \ref{theorem3} shows that all the higher order derivatives of $(g,F)$ are bounded uniformly. Thus they converge smoothly to some tuple $(g_T, F_T)$. Applying the short time existence of the flow starting at $(g_T,F_T)$, we see that the flow can be continued through $T$, contradicting the maximal existence time $T$.
\hfill \qed

%%%%%%%%%%%%%%%%%%%%%%%%%%%%%%%%%%%%%%%%

%%%%%%%%%%%%%%%%%%%%%%%%%%%%%%%%%%%%%%%%

\section{Estimates for the Dimensionally Reduced Flow}

%%%%%%%%%%%%%%%%%%%%%%%%%%

We come now to the proof of Theorem \ref{eqn:system}. Once again, part (a) is easy because $\p_td\beta= - d\Box_g\beta= - \Box_g d\beta$ and $\p_td\Psi= - d\Box_g\Psi= - \Box_g d\Psi$. Thus $d\beta$ and $d\Psi$ continue to vanish for all time if they vanish at time $t=0$. Part (b) is established in the same way as part (b) of Theorem \ref{Euclidean}: the leading terms in the flow for $f, \beta,\Psi$ are Laplacians with respect to the Riemannian metric $\hat g$, and hence strictly elliptic, while the leading term in the flow for $\hat  g_{ij}$ can also be made elliptic by DeTurck's trick as in Lemma \ref{lemma2}. The short-time existence of the flow follows again from the standard theory of parabolic partial differential equations.

\medskip

The main step is to establish part (c). For this, we
apply diffeomorphisms to simplify the evolution equations as follows: let $\varphi_t: M^{10-p}\to M^{10-p}$ be the one-parameter diffeomorphisms generated by $-\frac{p+1}{2}\nabla f$, and we consider the tuple $(g,f,\beta,\Psi)$ after pulled back by $\varphi_t$, e.g. $g_t = \varphi^*_t \hat g$, etc (abusing notations we still denote $(f,\beta,\Psi)$ to be the tuple after pulled back and use the notation $g_t$ for $\varphi^*_t\hat g$ since there is now no possibility of confusion with metrics on $M^{11}$). Using the formulas $\iota_{\nabla f} \beta = (-1)^p \star(df\wedge \star\beta)$ and $\iota_{\nabla f} \Psi = (-1)^p \star(df\wedge \star \Psi)$, we obtain
the following equations from (\ref{eqn:g1 main}), (\ref{eqn:f1 main}), (\ref{eqn:beta1 main}) and (\ref{eqn:psi1 main}), which are the main equations we will study in this section and in the equations the background metric for $\Delta$, $\nabla$, $|\cdot|^2$ or $\star$ is taken as the metric $g_t$:
\bea\label{eqn:g main}
\frac{\partial g_{ij}}{\partial t} &=& \nonumber- 2 R_{ij}+ \frac{p+1}{2} f_i f_j - \sigmaone e^{-(p+1)f} (\beta^2)_{ij} + (\Psi^2)_{ij} \\
&& \quad + \frac 1 3 (\sigmaone e^{-(p+1)f}|\beta|^2 - |\Psi|^2  )g_{ij}
\\
\label{eqn:f main}
\frac{\partial f}{\partial t} &=& \Delta f - \frac 2 3 \sigmaone e^{-(p+1)f} |\beta|^2 - \frac{1}{3}|\Psi|^2 - 2 \tilde \lambda e^{-f},\\
\label{eqn:beta main}
\frac{\partial\beta}{\partial t} &= & -\Box_g \beta + (-1)^{p+1}(p+1) d\star(df\wedge \star\beta)\nonumber \\
 &&\quad -\sigmaone  \frac{(-1)^{p+1} (p+1)}{4}  e^{\frac{p+1}{2}f} df\wedge \star(\Psi\wedge \Psi)  - \sigmaone \frac{(-1)^{p+1}}{2} e^{\frac{p+1}{2}f} d\star(\Psi\wedge \Psi),\\
\label{eqn:psi main}
\frac{\partial \Psi}{\partial t} &=& -\Box_g \Psi  -  \frac{p+1}{2} e^{-\frac{p+1}{2}f} df\wedge \star(\beta\wedge \Psi) + e^{-\frac{p+1}{2}f} d\star(\beta\wedge \Psi).
\eea
These equations are equivalent to the original ones, and it will suffice to derive the desired estimates for them. We assume the initial values are $(g_0,f_0,\beta_0,\Psi_0)$, where $f_0$ is a smooth function, $\beta_0$  a closed $(3-p)$-form and $\Psi_0$ a closed $4$-form on $M^{10-p}$. We may normalize the constant $\tilde \lambda$ to be $0$ or $\pm 1$ from now on.

%%%%%%%%%%%%%%%%%%%%

\medskip
Assume now that the flow exists on the time interval $[0,T)$, and that there exists a finite constant $K>0$, such that
\bea
\sup_{M^{10-p}\times [0,T)}(|Rm| + |f|+|\beta| + |\Psi|) \le K.
\eea
We are going to show that the higher order derivatives of $(g,f,\beta,\Psi)$ are bounded on the time interval $[0,T)$, by a constant depending only on $K, T$.
We start with some preliminary estimates.

\begin{lemma}\label{prop:1}
There is  a constant $C(K,T)>0$ such that
$$\sup_{M^{10-p}} |\nabla f|^2\le \frac{C}{t},\quad {\mathrm{for \,\, any\,\,}} t\in (0,T).$$
\end{lemma}
\noindent{\it Proof. }
By straightforward calculations, we have
\bea\nonumber
(\frac{\partial }{\partial t} - \Delta) |\nabla f|^2 &= & - 2 |\nabla^2 f|^2 - \frac{p+1}{2} |\nabla f|^4 + \sigmaone e^{-(p+1)f} (\beta^2)_{ij} f_i f_j - (\Psi^2)_{ij} f_i f_j\\
&& \nonumber +\sigmaone \frac{1}{6} e^{-(p+1) f} |\beta|^2 |\nabla f|^2 - \frac 16 |\Psi|^2 |\nabla f|^2 \\
&& \nonumber  - 2 \<  \nabla f, \nabla  ( - \sigmaone \frac 23 e^{-(p+1)f} |\beta|^2 - \frac 1 3 |\Psi|^2 - 2 \tilde \lambda e^{-f}   )  \>.
\eea
Thus for some constant $C_{K}>0$
\bea\label{eqn:nablaf2}
(\frac{\partial }{\partial t} - \Delta)|\nabla f|^2 \le - 2 |\nabla^2 f|^2 - \frac{p+1}{2}|\nabla f|^4 +C_{K}( |\nabla f|^2+ |\nabla f||\nabla \beta| + |\nabla f| |\nabla \Psi|  ) .
\eea
To control the terms involving $|\nabla \beta|$ and $|\nabla \Psi|$, we  calculate the equations for $|\beta|^2$ and $|\Psi|^2$. From (\ref{eqn:beta main}), we have
\bea\nonumber
&& ( \frac{\partial}{\partial t} - \Delta )|\beta|^2 \\ \nonumber = & & 2\Big( R_{ij} - \frac{p+1}{4} f_i f_j + \sigmaone \frac{1}{2} e^{-(p+1)f} (\beta^2)_{ij} - \frac{1}{2} (\Psi^2)_{ij} - \frac{1}{6}( \sigmaone e^{-(p+1)f} |\beta|^2 - |\Psi|^2  )g_{ij}  \Big) (\beta^2)_{ij}\\
&& \nonumber - 2 |\nabla \beta|^2 + Rm* \beta * \beta + \nabla^2 f * \beta * \beta + \nabla f * \nabla \beta + e^{\frac{p+1}{2}f} (\nabla f * \beta*\Psi * \Psi + \beta*\nabla \Psi * \Psi    )  %\\
%&& \nonumber+ e^{\frac{p+1}{2}f} (\nabla f * \beta*\Psi * \Psi + \beta*\nabla \Psi * \Psi    ),
\eea
so
\begin{equation}\label{eqn:beta2}(\frac{\partial}{\partial t} - \Delta)|\beta|^2 \le -2 |\nabla \beta|^2 + C_{K} ( |\nabla^2 f| + |\nabla f||\nabla \beta| + |\nabla f| + |\nabla \Psi| + 1  ).\end{equation}
Similarly we have
\bea\nonumber
\frac{\partial}{\partial t}|\Psi|^2  & =  & \Delta|\Psi|^2 - 2 |\nabla \Psi|^2 + 2 \Big( R_{ij} - \frac{p+1}{4} f_i f_j + \sigmaone \frac{1}{2} e^{-(p+1)f} (\beta^2)_{ij} - \frac{1}{2} (\Psi^2)_{ij} \\
&& \nonumber- \frac{1}{6}( \sigmaone e^{-(p+1)f} |\beta|^2 - |\Psi|^2  )g_{ij}  \Big) (\Psi^2)_{ij} + Rm* \Psi * \Psi + e^{-\frac{p+1}{2} f} \Big( \nabla f* \beta * \Psi * \Psi \\
&& \nonumber + \nabla \beta * \Psi *\Psi + \beta * \nabla \Psi * \Psi    \Big)\\
&\le & \Delta|\Psi|^2 - 2 |\nabla \Psi|^2  + C_{K} ( |\nabla f| + |\nabla \beta| + |\nabla \Psi| + 1   ).\label{eqn:psi2}
\eea
and
\begin{equation}\label{eqn:f sq}
(\frac{\partial}{\partial t} - \Delta) f^2 \le - 2 |\nabla f|^2 + C_{K}.
\end{equation}
Combining the equations (\ref{eqn:nablaf2}), (\ref{eqn:beta2}), (\ref{eqn:psi2}) and (\ref{eqn:f sq}), we get %\bea\nonumber
\begin{equation}\label{eqn:g0}(\frac{\partial}{\partial t} - \Delta) \Big( |\nabla f|^2 + |\beta|^2 + |\Psi|^2 \Big)  \le   -|\nabla^2 f|^2 - \frac{p+1}{2}|\nabla f|^4 - |\nabla \Psi|^2 - |\nabla \beta|^2 + C_{K},\end{equation}
and for $G_0:= t( |\nabla f|^2 + |\beta|^2 +|\Psi|^2  ) + A_0 f^2$
$$(\frac{\partial}{\partial t} - \Delta) G_0  \le   t \Big( -|\nabla^2 f|^2 - \frac{p+1}{2}|\nabla f|^4 - |\nabla \Psi|^2 - |\nabla \beta|^2 \Big) - |\nabla f|^2 + C_{K,T},$$
for $A_0>1$ sufficiently large depending only on $K, \, T$. In deriving the inequality above, we use the Cauchy-Schwarz inequalities, e.g. $|\nabla f| |\nabla \beta|\le \frac{1}{10}|\nabla \beta|^2 + \frac 52 |\nabla f|^2$. By maximum principle, it follows that $\sup_{M^{10-p}\times [0,T)} G_0\le C_{K,T}$, and from the definition of $G_0$, it follows that $\sup_{M^{10-p}} |\nabla f|^2 \le \frac{C_{K,T}}{t}$ for all $t\in (0,T)$.
%\eea
\hfill\qed

\begin{lemma}
The following formula holds for the Levi-Civita connection $\Gamma$
\begin{equation}\label{eqn:LC}
\frac{\partial}{\partial t} \Gamma = \nabla Rm + \nabla^2  f* \nabla f + e^{-(p+1)f} ( \nabla f* \beta * \beta + \nabla \beta * \beta   ) + \nabla \Psi * \Psi.
\end{equation}
\end{lemma}
\noindent{\it Proof. } We write the equation $\frac{\partial g_{ij}}{\partial t} = - 2 h_{ij}$, where  \begin{equation}\label{eqn:h ij}h_{ij}:= R_{ij} - \frac{p+1}{4} f_i f_j + \sigmaone \frac{1}{2} e^{-(p+1)f} (\beta^2)_{ij} - \frac{1}{2} (\Psi^2)_{ij}- \frac{1}{6}(\sigmaone e^{-(p+1)f} |\beta|^2 - |\Psi|^2  )g_{ij}.  \end{equation}
By definition $\Gamma^k_{ij} = \frac{1}{2} g^{kl}( \partial_i g_{jl} + \partial_j g_{il} - \partial_l g_{ij}   )$, it follows that
$$\frac{\partial }{\partial t}\Gamma^k_{ij} = - g^{kl}( \nabla_i h_{jl} + \nabla_j h_{il} - \nabla_l h_{ij}   ),$$
from which we get the desired equation for $\frac{\partial \Gamma}{\partial t}$. \hfill\qed

\begin{lemma}
The following equation holds for some $C_{K,T}>0$:
\bea\label{eqn:grad beta}
(\frac{\partial }{\partial t} - \Delta) |\nabla \beta|^2  & \le &  - 2 |\nabla^2 \beta|^2 + C_{K,T}\Big( |\nabla \beta|^2 + |\nabla Rm||\nabla \beta| + |\nabla^3 f| |\nabla f| \\ && \nonumber + |\nabla^2 f| |\nabla \beta|^2
 \nonumber+ |\nabla f| |\nabla\beta| |\nabla^2 \beta| + |\nabla f|^2 |\nabla \beta| + |\nabla^2 f| |\nabla \beta|\\
 & & \nonumber + |\nabla f| |\nabla \beta| |\nabla \Psi| + |\nabla^2 \Psi | |\nabla \beta| + |\nabla \Psi|^2 |\nabla \beta| +1   \Big)
\eea
\end{lemma}
\noindent{\it Proof. }
From the formula $\frac{\partial }{\partial t}\nabla \beta = \nabla \frac{\partial \beta}{\partial t} + \beta* \frac{\partial}{\partial t}\Gamma$, and the formulas  (\ref{eqn:LC}) and (\ref{eqn:beta main})
\bea\label{eqn:grad beta 1}
\frac{\partial}{\partial t}\nabla\beta & = & \Delta \nabla \beta + \nabla \beta * Rm + \beta * \nabla Rm + \nabla^3 f * \beta + \nabla^2 f * \nabla \beta \\
&& \nonumber + \nabla f * \nabla^2 \beta + e^{\frac{p+1}{2}f} \Big( \nabla f* \nabla f * \Psi * \Psi +  \nabla^2 f * \Psi * \Psi + \nabla f * \nabla \Psi * \Psi\\
&& \nonumber + \nabla^2 \Psi * \Psi + \nabla \Psi* \nabla \Psi\Big).
\eea
where we use the formula $\nabla \Delta \beta= \Delta \nabla \beta + \beta * \nabla Rm + \nabla \beta * Rm$. From this equation the estimate (\ref{eqn:grad beta}) follows.
\hfill \qed

\begin{lemma}
The following holds:
\bea\label{eqn:grad psi}
 (\frac{\partial }{\partial t} - \Delta)|\nabla \Psi|^2 & \le & - 2|\nabla^2 \Psi| ^2 + C_{K,T}\Big ( |\nabla \Psi|^2 + |\nabla Rm| |\nabla \Psi| + |\nabla \Psi|^2 \\ && \nonumber + |\nabla f|^2 |\nabla \Psi| | + |\nabla^2 f| |\nabla \Psi| + |\nabla f| |\nabla \beta| |\nabla \Psi| + |\nabla f| |\nabla \Psi|^2 \\
 && \nonumber+ |\nabla^2 \beta| |\nabla \Psi| + |\nabla \beta| |\nabla \Psi|^2 + |\nabla \Psi| |\nabla^2 \Psi| +1       \Big)
\eea
\end{lemma}
\noindent{\it Proof. } From the equation $\frac{\partial}{\partial t}\nabla \Psi = \nabla \frac{\partial}{\partial t}\Psi + \Psi * \frac{\partial}{\partial t}\Gamma$ and equations (\ref{eqn:psi main}), (\ref{eqn:LC}), we have
\bea\label{eqn:grad psi 1}
\frac{\partial \nabla \Psi}{\partial t} & = & \Delta \nabla \Psi + \nabla Rm * \Psi + Rm * \nabla \Psi + e^{-\frac{p+1}{2}f} \Big( \nabla f* \nabla f * \beta * \Psi  \\
&& \nonumber + \nabla^2 f* \beta * \Psi + \nabla f* \nabla \beta * \Psi + \nabla f * \beta * \nabla \Psi + \nabla^2 \beta * \Psi + \nabla \beta * \nabla \Psi + \beta * \nabla^2 \Psi  \Big).
\eea
The estimate (\ref{eqn:grad psi}) follows from the above equation.
\hfill\qed

\begin{lemma}
 $|Rm|^2$ satisfies the following
\begin{equation}\label{eqn:Rm 2}
(\frac{\partial }{\partial t} - \Delta)|Rm|^2 \le - 2|\nabla Rm|^2 + C_{K,T}\Big( |\nabla^2 f|^2 + |\nabla f|^2  + |\nabla \beta|^2 + |\nabla^2 \beta| + |\nabla\Psi|^2 + |\nabla^2 \Psi| + 1     \Big)
\end{equation}
\end{lemma}
\noindent{\it Proof. } From the evolution equation of $Rm$ that
%\bea\nonumber
$$\frac{\partial}{\partial t} R_{ijkl} = - \nabla_i \nabla_l h_{jk} - \nabla_j \nabla_k h_{il} + \nabla_i \nabla_k h_{jl} + \nabla_j \nabla_l h_{ik} - h_{km} R_{ijml} - h_{ml}R_{ijkm} ,$$
where $h_{ij}$ is given in (\ref{eqn:h ij}), we get
\bea \nonumber
\frac{\partial Rm}{\partial t} & = & \Delta Rm + Rm * Rm + \nabla^2 f* \nabla^2 f + Rm * \Psi * \Psi + \nabla \Psi * \nabla \Psi + \nabla^2 \Psi * \Psi\\
&& \nonumber + e^{-(p+1)f } \Big( Rm* \beta* \beta + \nabla f * \nabla f * \beta * \beta + \nabla^2 f* \beta * \beta + \nabla f * \nabla \beta * \beta \\
&& + \nabla \beta * \nabla \beta + \nabla^2 \beta * \beta    \Big), \label{eqn:Rm eqn}
\eea
where we apply Riccati equations in deriving the above, in particular,  $\nabla^3 f$ terms can be cancelled in the equations. The estimate (\ref{eqn:Rm 2}) follows from  (\ref{eqn:Rm eqn}). \hfill \qed

\begin{lemma}
We have the following equation for $|\nabla^2 f|^2$
\bea\label{eqn:nabla2f}
(\frac{\partial}{\partial t} - \Delta) |\nabla^2 f|^2% & \le & -2 |\nabla^3 f|^2 + C_{K,T}\Big( |\nabla^2 f|^2 + |\nabla^2 f| |\nabla^2 \beta| + |\nabla^2 f| |\nabla \beta|^2 \\
%&& \nonumber + |\nabla^2 f| |\nabla f| |\nabla \beta| + |\nabla ^2 f | |\nabla^2 \Psi| + |\nabla^2 f| |\nabla \Psi|^2 + |\nabla^2 f| |\nabla f| |\nabla \Psi| + 1   \Big)\\
& \le & \nonumber - 2 |\nabla^3 f|^2 + C_{K,T} \Big( |\nabla^2 f|^2 + |\nabla f|^4 + 1   \Big)\\
&& \nonumber +\frac{1}{10} |\nabla^2 \beta|^2 + \frac{1}{10}|\nabla \beta|^4 + \frac{1}{10}|\nabla^2 \Psi|^2 + \frac{1}{10}|\nabla\Psi|^4
\eea
\end{lemma}
\noindent{\it Proof. } From the equation $f_{ij} = \partial^2_{ij} f - \Gamma_{ij}^k f_k$ and the evolution equation (\ref{eqn:f main}) for $f$, we have
\bea\nonumber
\frac{\partial}{\partial t} \nabla^2 f & = & \nabla^2 f + \nabla f * \frac{\partial \Gamma}{\partial t}\\
& = &\nonumber \Delta \nabla^2 f + \nabla^2 f * Rm + \nabla^2 \Psi * \Psi + \nabla\Psi * \nabla \Psi - \frac{p+1}{2}|\nabla f|^2 \nabla^2 f\\
&& \nonumber + \nabla f * \nabla \Psi * \Psi  + \tilde \lambda e^{-f} ( \nabla^2 f + \nabla f* \nabla f  ) \\
&& \nonumber + e^{-(p+1)f} \Big( \nabla^2 \beta* \beta + \nabla \beta * \nabla \beta + \nabla f* \beta * \nabla \beta + \nabla^2 f * \beta * \beta \\
&&+ \nabla f * \nabla f * \beta * \beta   \Big), \label{eqn:nabla2f 1}
\eea
from (\ref{eqn:nabla2f 1}) we have
\bea\nonumber%\label{eqn:nabla2f}
(\frac{\partial}{\partial t} - \Delta) |\nabla^2 f|^2  & \le & -2 |\nabla^3 f|^2 + C_{K,T}\Big( |\nabla^2 f|^2 + |\nabla^2 f| |\nabla^2 \beta| + |\nabla^2 f| |\nabla \beta|^2 \\
&& \nonumber + |\nabla^2 f| |\nabla f| |\nabla \beta| + |\nabla ^2 f | |\nabla^2 \Psi| + |\nabla^2 f| |\nabla \Psi|^2 + |\nabla^2 f| |\nabla f| |\nabla \Psi| + 1   \Big)\\
& \le & \nonumber - 2 |\nabla^3 f|^2 + C_{K,T} \Big( |\nabla^2 f|^2 + |\nabla f|^4 + 1   \Big) + \frac{1}{10}|\nabla^2 \beta|^2 \\
&& \nonumber + \frac{1}{10}|\nabla \beta|^4 + \frac{1}{10}|\nabla^2 \Psi|^2 + \frac{1}{10}|\nabla\Psi|^4
\eea
by Cauchy-Schwarz inequalities. \hfill\qed

\medskip

With these formulas in hand, we are now ready to prove the derivatives estimates.
\begin{lemma}\label{prop:2}
There exists a constant $C(K,T)>0$ such that
\begin{equation}\label{eqn:1st estimate}
\sup_{M^{10-p}}\Big(  |\nabla\beta |^2 + |\nabla \Psi|^2 + |\nabla^2 f|^2\Big) \le \frac{C}{t},\quad \forall ~ t\in (0,T).
\end{equation}
\end{lemma}
\noindent{\it Proof. } By (\ref{eqn:beta2}) and (\ref{eqn:grad beta}), we have for any constant $A_2>0$
\bea\label{eqn:prop2 1}
&& (\frac{\partial}{\partial t} - \Delta) \Big( (|\beta|^2 + A_2) |\nabla \beta|^2\Big)\\
&= & \nonumber (\frac{\partial}{\partial t} - \Delta)|\beta|^2 |\nabla \beta|^2 + (|\beta|^2+A_2) (\frac{\partial}{\partial t} - \Delta) |\nabla \beta|^2 - 2 \< \nabla |\beta|^2 ,\nabla |\nabla \beta|^2 \>\\
& \le &\nonumber -2 |\nabla \beta|^4  - 2 A_2 |\nabla^2 \beta|^2 + C_{K,T}\Big(1+ |\nabla \beta|^2 + |\nabla \beta|^2 |\nabla^2 f| + |\nabla f| |\nabla \beta|^3 + |\nabla \Psi | |\nabla \beta|^2 \\
&& \nonumber + A_2 |\nabla Rm| |\nabla \beta| + A_2 |\nabla^3 f| |\nabla \beta| + A_2 |\nabla^2 f| |\nabla \beta|^2 + A_2 |\nabla f| |\nabla \beta| |\nabla^2 \beta| + A_2 |\nabla f|^2 |\nabla \beta| \\
&& \nonumber+ A_2 |\nabla^2 f| |\nabla \beta| + A_2 |\nabla f| |\nabla \beta| |\nabla \Psi| + A_2 |\nabla^2 \Psi| |\nabla \beta| + A_2 |\nabla \Psi|^2 |\nabla \beta|  + |\nabla\beta|^2 |\nabla^2 \beta|  \Big)\\
&\le & \nonumber- \frac 3 2 |\nabla \beta|^4 - A_2 |\nabla^2 \beta|^2 + C_{K,T}\Big( 1+ |\nabla^2 f|^2 + |\nabla f|^4  \Big) + \frac{1}{10}|\nabla Rm|^2     \\
&&\nonumber+ \frac{1}{10} |\nabla ^3 f|^2  + \frac{1}{10}|\nabla^2 \Psi|^2 + \frac{1}{10}|\nabla \Psi|^4,
\eea
for $A_2 >1$ large enough depending only on $K,\, T$, and in the last inequality above, we use Cauchy-Schwarz inequalities to simplify the expression.

By (\ref{eqn:psi2}) and (\ref{eqn:grad psi}), we have for any constant $A_1>0$
\bea\label{eqn:prop2 2}
&& (\frac{\partial}{\partial t} - \Delta) \Big( (|\Psi|^2 + A_1) |\nabla \Psi|^2  \Big)\\
& \le  & \nonumber- 2 |\nabla \Psi|^4 - 2A_1 |\nabla^2 \Psi|^2 + C_{K,T} \Big( 1+ A_1 |\nabla\Psi|^2 + A_1 |\nabla f| |\nabla \Psi|^2 +A_1 |\nabla Rm| |\nabla \Psi| \\
&& \nonumber + A_1 |\nabla f|^2 |\nabla \Psi| + A_1 |\nabla^2 f| |\nabla \Psi|  + A_1 |\nabla f | |\nabla \beta| |\nabla \Psi|  + A_1 |\nabla^2 \beta |  |\nabla \Psi| \\
&& \nonumber + A_1 |\nabla \beta| |\nabla \Psi|^2 + |\nabla \Psi| |\nabla^2 \Psi|     \Big)\\
&\le & \nonumber - |\nabla\Psi|^4 - A_1 |\nabla^2 \Psi|^2 + C_{K,T}\Big(1  + |\nabla f|^4 + |\nabla^2 f|^2 + |\nabla \beta|^2  \Big) + \frac{1}{10}|\nabla^2 \beta|^2 + \frac{1}{10}|\nabla Rm|^2,
\eea
where in the last inequality we take $A_1>1$ to be sufficiently large and apply the Cauchy-Schwarz inequalities. We denote
\begin{equation}\label{eqn:G1}G_1: =  (|\Psi|^2 + A_1) |\nabla\Psi|^2 + (|\beta|^2 + A_2)|\nabla \beta|^2 + |Rm|^2 + |\nabla ^2 f|^2.\end{equation}
By (\ref{eqn:Rm 2}),  (\ref{eqn:nabla2f}), (\ref{eqn:prop2 1}) and (\ref{eqn:prop2 2}), after some cancellations  we have
\bea\nonumber
&&(\frac{\partial}{\partial t} - \Delta ) G_1\\
&\le & \nonumber -  |\nabla \beta|^4 - |\nabla^2 \beta|^2 - |\nabla \Psi|^4 - |\nabla^2 \Psi|^2 - |\nabla Rm|^2 - |\nabla^3 f|^2 \\
&& \nonumber+ C_{K,T}\Big(1 + |\nabla^2 f|^2 +|\nabla f|^4   \Big).
\eea
Combining (\ref{eqn:g0}) we get
\bea\nonumber
&&(\frac{\partial}{\partial t} - \Delta) \Big( t G_1 + A_3 (|\nabla f|^2 +|\beta|^2 + |\Psi|^2) \Big)\\
& \le &\nonumber   G_1 + t \Big(-  |\nabla \beta|^4 - |\nabla^2 \beta|^2 - |\nabla \Psi|^4 - |\nabla^2 \Psi|^2 - |\nabla Rm|^2 - |\nabla^3 f|^2 \\
&& \nonumber+ C_{K,T}(1 + |\nabla^2 f|^2 +|\nabla f|^4 )\Big) - A_3  ( |\nabla^2 f|^2 + |\nabla f|^4 + |\nabla \Psi|^2 + |\nabla \beta|^2   )+ A_3C_{K,T}\\
&\le & \label{eqn:need1}  t \Big( - |\nabla^2 \beta|^2 - |\nabla^2 \Psi|^2 - |\nabla Rm|^2 - |\nabla^3 f|^3  \Big) + C_{K,T} A_3 + C_{K,T}\\
&\le & \nonumber C_{K,T}+ A_3 C_{K,T},
\eea
if we choose $A_3>1$ large enough which depends only on $K,\, T$. Applying maximum principle, we get $\sup_{M^{10-p}} G_1 \le \frac{C_{K,T}}{t}$, and hence
$$\sup_{M^{10-p}} \Big( |\nabla \beta|^2 + |\nabla \Psi|^2 + |\nabla^2 f|^2  \Big)\le \frac{C_{K,T}}{t},\quad \forall~ t\in (0,T).$$
This finishes the proof.
\hfill \qed

\begin{lemma}
$|\nabla Rm|^2$ satisfies the following inequality
\bea\label{eqn:grad Rm}
&& \ho |\nabla Rm|^2\\
& \le &\nonumber - 2 |\nabla^2 Rm|^2+ C_{K,T} |\nabla Rm| \Big( |\nabla Rm| + t^{-1/2}( |\nabla ^3 f| + |\nabla^2 \beta| + |\nabla^2 \Psi|  ) \\
&& \quad + t^{-1/2} + t^{-3/2} + t^{-1} + |\nabla^3 f|
 \nonumber  + |\nabla^3 \beta|  + |\nabla^3\Psi|  \Big).
\eea
\end{lemma}
\noindent{\it Proof. } We calculate using the formula $\frac{\partial}{\partial t}\nabla Rm = \nabla \frac{\partial Rm}{\partial t} + Rm * \frac{\partial\Gamma}{\partial t}$ and (\ref{eqn:Rm eqn})
\bea\label{eqn:grad Rm 1}
&& \frac{\partial}{\partial t} \nabla Rm \\ &=&\nonumber \Delta \nabla Rm + Rm * \nabla Rm+ \nabla^2 f * \nabla^3 f + \nabla Rm * \Psi *\Psi + Rm* \nabla \Psi * \Psi + Rm * \nabla^2 f * \nabla f\\
&&\nonumber + \nabla^2 \Psi * \nabla \Psi + \nabla^3 \Psi * \Psi + e^{-(p+1)f}\Big( Rm * \nabla f * \beta * \beta + \nabla Rm * \beta*\beta + Rm * \nabla \beta * \beta \\
&& \nonumber + \nabla f * \nabla f * \nabla f * \beta * \beta + \nabla^2 f * \nabla f * \beta * \beta +  \nabla f * \nabla f * \nabla \beta * \beta + \nabla^3 f * \beta * \beta \\
&&\nonumber + \nabla^2 f * \nabla \beta * \beta+ \nabla f* \nabla^2 \beta * \beta + \nabla f * \nabla \beta * \nabla \beta + \nabla^2 \beta * \nabla \beta + \nabla^3 \beta * \beta   \Big),
\eea
therefore it follows that
\bea\nonumber
&& \ho |\nabla Rm|^2\\
& \le &\nonumber - 2 |\nabla^2 Rm|^2 + C_{K,T} |\nabla Rm |\Big( |\nabla Rm| + |\nabla^2 f| |\nabla^3 f| + |\nabla f| + |\nabla \beta| + |\nabla \Psi|  \\
&& \nonumber \quad+ |\nabla f|^3 + |\nabla f| |\nabla ^2 f| + |\nabla f|^2 |\nabla \beta| + |\nabla^3 f| + |\nabla^2 f| |\nabla \beta| + |\nabla f| |\nabla^2 \beta| \\
&& \nonumber \quad + |\nabla f| |\nabla \beta| + |\nabla \beta| |\nabla^2 \beta| + |\nabla^3 \beta| + |\nabla^2 \Psi| |\nabla \Psi| + |\nabla^3 \Psi|  \Big)\\
&\le &\nonumber - 2 |\nabla^2 Rm|^2+ C_{K,T} |\nabla Rm| \Big( |\nabla Rm| + t^{-1/2} |\nabla ^3 f| + t^{-1/2} + t^{-3/2} + t^{-1} + |\nabla^3 f| \\
& & \nonumber \quad  + t^{-1/2} |\nabla^2 \beta| + |\nabla^3 \beta| + t^{-1/2} |\nabla^2 \Psi| + |\nabla^3\Psi|  \Big).
\eea
\hfill \qed
\begin{lemma}\label{Lemma for nabla3f}
$|\nabla^3 f|^2$ satisfies the inequality
\bea\label{eqn:nabla3f}
&& \ho |\nabla^3 f|^2\\
&\le & \nonumber - 2 |\nabla^4 f|^2 + C_{K,T}|\nabla^3 f| \Big( |\nabla^3 f| + t^{-1/2} (|\nabla Rm| + |\nabla^2 \beta| + |\nabla^2 \Psi| )+ t^{-1} + t^{-3/2} \\
&&\nonumber \quad  +|\nabla^3 \beta| +|\nabla^3 \Psi| +  + t^{-1} |\nabla^3 f|  \Big).
\eea
\end{lemma}
\noindent{\it Proof. }  By (\ref{eqn:nabla2f 1}) we calculate
\bea\nonumber
\frac{\partial }{\partial t} \nabla^3 f & =   & \nabla \frac{\partial}{\partial t}\nabla^2 f + \nabla^2 f* \frac{\partial \Gamma}{\partial t}\\
& = &\nonumber\Delta \nabla^3 f + \nabla^3 f * Rm + \nabla^2 f * \nabla Rm + e^{-(p+1)f}\Big( \nabla^3 f * \beta * \beta + \nabla^2 f *\nabla \beta * \beta\\
 && \nonumber  + \nabla f * \nabla f * \nabla \beta * \beta + \nabla f * \nabla^2 \beta * \beta + \nabla f * \nabla \beta * \nabla \beta + \nabla^3 \beta * \beta + \nabla^2 \beta* \nabla \beta\\
 && \nonumber + \nabla f * \nabla f* \nabla f*\beta * \beta + \nabla^2 f * \nabla f * \beta * \beta + \nabla^2 f * \nabla \beta * \beta \Big) + \nabla^3 \Psi * \Psi\\
 &&\nonumber + \nabla^2 \Psi * \nabla\Psi + \nabla^2 f * \nabla^2 f * \nabla f + \nabla^3 f * \nabla f*\nabla f + \nabla^2 f * \nabla \Psi * \Psi\\
 && \nonumber + \nabla f * \nabla^2 \Psi * \Psi + \nabla f * \nabla \Psi *\nabla\Psi + \tilde \lambda  e^{-f}( \nabla^2 f * \nabla f + \nabla^3 f + \nabla f*\nabla f* \nabla f   ),
\eea
Thus by the Cauchy-Schwarz inequality and Lemmas \ref{prop:1} and \ref{prop:2}, we have
\bea\nonumber
&& \ho |\nabla^3 f|^2 \\
& \le & \nonumber - 2 |\nabla^4 f|^2 + C_{K,T}|\nabla ^3 f| \Big( |\nabla^3 f| + |\nabla^2 f||\nabla Rm| + |\nabla^2 f| |\nabla \beta| + |\nabla f|^2 |\nabla\beta| + |\nabla f| |\nabla ^2 \beta| \\
&&\nonumber \quad + |\nabla f| |\nabla\beta|^2 + |\nabla^3\beta| + |\nabla^2 \beta| |\nabla \beta| + |\nabla f|^3 + |\nabla^3 \Psi| + |\nabla \Psi| |\nabla^2 \Psi| + |\nabla^2 f | |\nabla f| \\
&&\nonumber \quad+ |\nabla f| |\nabla^2 f|^2 + |\nabla f|^2 |\nabla^3 f| + |\nabla^2 f| |\nabla \Psi| + |\nabla f| |\nabla^2 \Psi| + |\nabla f| |\nabla \Psi|^2   \Big)\\
&\le &\nonumber - 2 |\nabla^4 f|^2 + C_{K,T}|\nabla^3 f| \Big( |\nabla^3 f| + t^{-1/2} (|\nabla Rm| + |\nabla^2 \beta| + |\nabla^2 \Psi| )+ t^{-1} + t^{-3/2} \\
&&\nonumber \quad  +|\nabla^3 \beta| +|\nabla^3 \Psi| +  + t^{-1} |\nabla^3 f|  \Big).
\eea
This finishes the proof of Lemma \ref{Lemma for nabla3f}.
\hfill\qed

\begin{lemma}Along the flow
$|\nabla^2\beta|^2$ satisfies the inequality
\bea\label{eqn:nabla2beta}
&& \ho|\nabla^2 \beta|^2\\
 &\le & \nonumber -2 |\nabla^3 \beta|^2 + C_{K,T}|\nabla^2 \beta|\Big( |\nabla^2 \beta| + t^{-1/2} ( |\nabla Rm| + |\nabla^3 f| + |\nabla^2 \beta| + |\nabla^2 \Psi| ) + |\nabla^2 Rm| \\
&& \nonumber\quad + |\nabla^4 f|
+ t^{-1/2} |\nabla^3 \beta| + t^{-1} |\nabla^2 \beta| + t^{-1} + t^{-3/2}  + |\nabla^3 f| + |\nabla^3 \Psi|  \Big),
\eea for some constant $C_{K,T}>0$.
\end{lemma}
\noindent{\it Proof. } By (\ref{eqn:grad beta 1}), we have
\bea\nonumber
\frac{\partial\nabla^2 \beta}{\partial t} & = & \nabla \frac{\partial\nabla \beta}{\partial} + \nabla \beta * \frac{\partial\Gamma}{\partial t}\\
& = &\nonumber \Delta \nabla^2 \beta  + \nabla^2 \beta * Rm + \nabla \beta * \nabla Rm + \beta* \nabla^2 Rm + \nabla^4 f * \beta + \nabla^3 f * \nabla \beta\\
&& \nonumber \quad + \nabla^2 f * \nabla^2 \beta + \nabla f* \nabla^3 \beta + \nabla^2 f*\nabla f * \nabla \beta + \nabla \beta * \nabla \Psi * \Psi\\
&& \nonumber \quad + e^{\frac{p+1}{2} f} \Big( \nabla f* \nabla f * \nabla f* \Psi * \Psi + \nabla f* \nabla^2 f * \Psi * \Psi + \nabla f* \nabla f * \nabla \Psi * \Psi \\
&& \nonumber\quad  + \nabla f*\nabla^2 \Psi * \Psi + \nabla f*\nabla \Psi *\nabla \Psi + \nabla^3 f * \Psi * \Psi + \nabla^2 f * \nabla \Psi * \Psi \\
&&\nonumber\quad  + \nabla^3 \Psi * \Psi + \nabla^2 \Psi * \nabla\Psi   \Big) + e^{-(p+1)f} ( \nabla f * \nabla \beta * \beta * \beta + \nabla \beta * \nabla \beta * \beta  ).
\eea
Then by Cauchy-Schwarz inequality and Lemmas \ref{prop:1} and   \ref{prop:2}, we have
\bea\nonumber
&& \ho |\nabla^2 \beta|^2 \\
& \le & \nonumber -2 |\nabla^3\beta|^2 + C_{K,T} |\nabla^2 \beta| \Big( |\nabla^2 \beta| + |\nabla \beta| |\nabla Rm| + |\nabla^2 Rm| + |\nabla^4 f| + |\nabla^3 f| |\nabla \beta| \\
&&\nonumber \quad+ |\nabla^2 \beta| |\nabla^2 f| + |\nabla f| |\nabla^3 \beta| +  |\nabla f| |\nabla \beta| |\nabla^2 \beta| + |\nabla\beta| |\nabla \Psi| + |\nabla f|^3 + |\nabla f| |\nabla^2 f| \\
&&\nonumber \quad+ |\nabla f|^2 |\nabla \Psi| + |\nabla f| |\nabla^2 \Psi| + |\nabla f| |\nabla \Psi|^2 + |\nabla^3 f| + |\nabla^2 f| |\nabla \Psi| + |\nabla^3 \Psi|\\
&&\nonumber \quad + |\nabla^2 \Psi| |\nabla \Psi| + |\nabla f| |\nabla\beta| + |\nabla \beta|^2   \Big)\\
& \le & \nonumber - 2 |\nabla^3 \beta|^2 + C_{K,T}|\nabla^2 \beta|\Big( |\nabla^2 \beta| + t^{-1/2} ( |\nabla Rm| + |\nabla^3 f| + |\nabla^2 \beta| + |\nabla^2 \Psi| ) + |\nabla^2 Rm| \\
&& \nonumber\quad + |\nabla^4 f|
+ t^{-1/2} |\nabla^3 \beta| + t^{-1} |\nabla^2 \beta| + t^{-1} + t^{-3/2}  + |\nabla^3 f| + |\nabla^3 \Psi|  \Big).
\eea
\hfill\qed
%%%new lemma
\begin{lemma}Along the flow
$|\nabla^2 \Psi|^2$ satisfies the inequality
\bea\label{eqn:nabla2psi}
&& \ho |\nabla^2 \Psi|^2 \\
& \le & \nonumber -2 |\nabla^3 \Psi|^2 + C_{K,T}|\nabla^2 \Psi|\Big( |\nabla^2 \Psi| + t^{-1/2} (|\nabla Rm| + |\nabla^2 \beta| + |\nabla^2 \Psi| )+ |\nabla^2 Rm|  \\
&& \nonumber\quad + t^{-3/2} + t^{-1}+|\nabla^3 f|  + |\nabla^3 \beta| + |\nabla^3 \Psi|    \Big),
\eea for some constant $C_{K,T}>0$.
\end{lemma}
\noindent{\it Proof. } By (\ref{eqn:grad psi 1}) we have the equation
\bea\nonumber
&& \frac{\partial}{\partial t} \nabla^2 \Psi\\
&= & \nonumber \nabla \frac{\partial\nabla\Psi}{\partial t } + \nabla \Psi * \frac{\partial \Gamma}{\partial t }\\
& = & \nonumber \Delta \nabla^2 \Psi + \nabla^2 \Psi * Rm + \nabla\Psi * \nabla Rm + \nabla^2 Rm * \Psi + \nabla \Psi * \nabla^2 f * \nabla f \\
&&\nonumber + e^{-\frac{p+1}{2}f} \Big( \nabla f*\nabla f * \nabla f * \beta * \Psi + \nabla f * \nabla^2 * \beta * \Psi + \nabla f* \nabla f * \nabla \beta * \Psi \\
&&\nonumber \quad+ \nabla f* \nabla f * \beta * \nabla \Psi + \nabla f * \nabla^2 \beta * \Psi + \nabla f* \nabla \beta * \nabla \Psi + \nabla f * \beta * \nabla ^2 \Psi \\
&&\nonumber \quad + \nabla^3 f * \beta * \Psi + \nabla^2 f * \beta * \nabla \Psi + \nabla^3 \beta * \Psi + \nabla^2 \beta * \nabla \Psi + \nabla \beta * \nabla^2\Psi + \beta * \nabla^3\Psi   \Big)\\
&&\nonumber + e^{-(p+1)f} ( \nabla f * \nabla \Psi * \beta * \beta + \nabla \beta * \beta * \nabla \Psi + \nabla\Psi * \nabla\Psi * \Psi  ).
\eea
By Cauchy-Schwarz inequality, Lemmas  \ref{prop:1} and  \ref{prop:2} we have
\bea\nonumber
&& \ho |\nabla^2 \Psi|^2 \\
& \le & \nonumber - 2 |\nabla^3\Psi|^2 + C_{K,T} |\nabla^2 \Psi| \Big( |\nabla^2 \Psi| + |\nabla \Psi| |\nabla Rm| + |\nabla^2 Rm| + |\nabla \Psi| |\nabla f| |\nabla^2 f| + |\nabla^3 f|  \\
&&\nonumber \quad+ |\nabla f| |\nabla^2 f| + |\nabla f|^2 |\nabla \beta| + |\nabla f|^2 |\nabla \Psi| + |\nabla f| |\nabla^2 \beta| + |\nabla f| |\nabla \beta| |\nabla \Psi| \\
&&\nonumber \quad+ |\nabla f| |\nabla^2 \Psi| + |\nabla^3 f| + |\nabla^2 f| |\nabla \Psi| + |\nabla^3 \beta| + |\nabla^2 \beta| |\nabla \Psi|  + |\nabla \beta| |\nabla^2 \Psi| + |\nabla^3 \Psi| \\
&&\nonumber \quad + |\nabla f| |\nabla \Psi| + |\nabla \beta| |\nabla \Psi| + |\nabla \Psi|^2  \Big)\\
& \le &\nonumber -2 |\nabla^3 \Psi|^2 + C_{K,T}|\nabla^2 \Psi|\Big( |\nabla^2 \Psi| + t^{-1/2} (|\nabla Rm| + |\nabla^2 \beta| + |\nabla^2 \Psi| )+ |\nabla^2 Rm|  \\
&& \nonumber\quad + t^{-3/2} + t^{-1}+|\nabla^3 f|  + |\nabla^3 \beta| + |\nabla^3 \Psi|    \Big).
\eea
\hfill\qed
\begin{lemma}\label{prop:3}
There exists a constant $C = C(K,T)>0$ such that
$$\sup_{M^{10-p}} \Big( |\nabla Rm| + |\nabla^2 \beta|  + |\nabla^2 \Psi|+ |\nabla^3f|  \Big)\le \frac{C}{t},\quad \forall ~ t\in (0,T)$$
\end{lemma}
%\newpage
\noindent{\it Proof. } Combining the inequalities (\ref{eqn:grad Rm}), (\ref{eqn:nabla3f}), (\ref{eqn:nabla2beta}) and (\ref{eqn:nabla2psi}) and applying the Cauchy-Schwarz inequality serval times, it follows that for $G_2: = |\nabla Rm|^2 + |\nabla^3 f|^2 + |\nabla^2 \beta|^2 + |\nabla^2 \Psi|^2$,
\bea\nonumber
&& \ho (t^2 G_2 )\\
& \le &\nonumber 2 t G_2 + t^2 \Big\{ -|\nabla^2 Rm|^2 - |\nabla^4 f|^2 - |\nabla^3 \beta|^2 \\
&& \quad- |\nabla^3 \Psi|^2
 \nonumber  + C_{K,T} (1+t^{-1/2})\Big(  |\nabla Rm|^2 + |\nabla^3 f|^2 + |\nabla^2 \beta|^2 + |\nabla^2 \Psi|^2 \Big) \\
 && \nonumber \quad + C_{K,T} t^{-3/2} ( |\nabla Rm| + |\nabla^3 f| + |\nabla ^2 \beta| + |\nabla^2 \Psi|  ) + C_{K,T} t^{-1} |\nabla^2 \beta|^2 +C_{K,T}     \Big\} \\
& \le &\nonumber C_{K,T} t G_2+ t^2 \Big\{ -|\nabla^2 Rm|^2 - |\nabla^4 f|^2   - |\nabla^3 \beta|^2 - |\nabla^3 \Psi|^2
 \nonumber \Big\} + C_{K,T}. %(1+t^{-1/2})\Big(  |\nabla Rm|^2 + |\nabla^3 f|^2 + |\nabla^2 \beta|^2 + |\nabla^2 \Psi|^2 \Big) \\
% && \nonumber + C_{K,T} t^{-3/2} ( |\nabla Rm| + |\nabla^3 f| + |\nabla ^2 \beta| + |\nabla^2 \Psi|  ) + C_{K,T} t^{-1} |\nabla^2 \beta|^2 +C_{K,T}     \Big\}
\eea
Combining with the inequality (\ref{eqn:need1}) with $G_1$ given in (\ref{eqn:G1}), we get for $A_4$ large enough
\bea\nonumber
&& \ho \Big\{t^2 G_2 + A_4\Big( t G_1 + A_3 (|\nabla f|^2 +|\beta|^2 + |\Psi|^2)\Big) \Big\}\\
&\le & \nonumber  t^2 \Big( -|\nabla^2 Rm|^2 - |\nabla^4 f|^2 - |\nabla^3 \beta|^2 - |\nabla ^3 \Psi|^2  \Big) + C_{K,T}.
\eea
Applying maximum principle it follows that $\sup_{M^{10-p}} G_2 \le \frac{C_{K,T}}{t^2}$. The desired estimate then follows from the definition of $G_2$. \hfill \qed

%\textcolor{black}{
%\bea
%&&\nonumber + t^2\Big\{- 2 |\nabla^2 Rm|^2+ C_{K,T} |\nabla Rm| \Big( |\nabla Rm| + t^{-1/2}( |\nabla ^3 f| + |\nabla^2 \beta| + |\nabla^2 \Psi|  ) \\
%&& \quad + t^{-1/2} + t^{-3/2} + t^{-1} + |\nabla^3 f|
 %\nonumber  + |\nabla^3 \beta|  + |\nabla^3\Psi|  \Big)\\
%&& \nonumber - 2 |\nabla^4 f|^2 + C_{K,T}|\nabla^3 f| \Big( |\nabla^3 f| + t^{-1/2} (|\nabla Rm| + |\nabla^2 \beta| + |\nabla^2 \Psi| )+ t^{-1} + t^{-3/2} \\
%&&\nonumber \quad  +|\nabla^3 \beta| +|\nabla^3 \Psi| +  + t^{-1} |\nabla^3 f|  \Big )\\
%&& \nonumber -2 |\nabla^3 \beta|^2 + C_{K,T}|\nabla^2 \beta|\Big( |\nabla^2 \beta| + t^{-1/2} ( |\nabla Rm| + |\nabla^3 f| + |\nabla^2 \beta| + |\nabla^2 \Psi| ) + |\nabla^2 Rm| \\
%&& \nonumber\quad + |\nabla^4 f|
%+ t^{-1/2} |\nabla^3 \beta| + t^{-1} |\nabla^2 \beta| + t^{-1} + t^{-3/2}  + |\nabla^3 f| + |\nabla^3 \Psi|  \Big)\\
%&&\nonumber -2 |\nabla^3 \Psi|^2 + C_{K,T}|\nabla^2 \Psi|\Big( |\nabla^2 \Psi| + t^{-1/2} (|\nabla Rm| + |\nabla^2 \beta| + |\nabla^2 \Psi| )+ |\nabla^2 Rm|  \\
%&& \nonumber\quad + t^{-3/2} + t^{-1}+|\nabla^3 f|  + |\nabla^3 \beta| + |\nabla^3 \Psi|    \Big)\Big\}
%\eea}

Next we prove the higher order estimates for general $m\in {\mathbf Z}_+$.

\begin{theorem}\label{prop:4}
For any $2\le m\in {\mathbf Z}_+$, there exists a constant $C=C(m,K,T)>0$ such that
\begin{equation}\label{eqn:grad m}\sup_{M^{10-p}} \Big( |\nabla^{m-1} Rm| + |\nabla^m \beta| + |\nabla^m \Psi| + |\nabla^{m+1} f|    \Big)\le \frac{C(m,K,T)}{t^{m/2}},\quad \forall~ t\in (0,T).\end{equation}
\end{theorem}

We will prove this theorem by induction. To begin with, we need equations on the higher order derivatives of $Rm$, $f$, $\beta$ and $\Psi$. For notation convenience we denote $$\gamma = \nabla f* \nabla f +  e^{-(p+1)f} \beta * \beta + \Psi* \Psi,$$ by (\ref{eqn:LC}), we know $\frac{\partial \Gamma}{\partial t} = \nabla Rm + \nabla \gamma$.

\begin{lemma}\label{lemma beta}
We have the following evolution equation for $\nabla^m \beta$
\bea\label{eqn:grad beta m}
\frac{\partial\nabla^m \beta}{\partial t} & =  & \Delta \nabla^m\beta + \sum_{i+j = m} \nabla^i \beta * \nabla^ j Rm + \sum_{i+j=m+1} \nabla^{i+1} f * \nabla^j \beta\\
&& \nonumber + \sum_{i+j+k = m+1} \nabla^i (e^{\frac{p+1}{2}f}) * \nabla^j\Psi * \nabla^k \Psi + \sum_{i+j = m-1}\nabla^i \beta * \nabla^{j+1}\gamma.
\eea
where in the summations $i, j, k$ are nonnegative integers.
\end{lemma}
\noindent{\it Proof. } We will prove this formula by induction. The formula holds for $m=1$ by (\ref{eqn:grad beta 1}). Assume it has been proved for $m-1$, then
\bea\nonumber
&&\frac{\partial\nabla^m \beta}{\partial t}\\
\nonumber & = & \nabla \frac{\partial\nabla^{m-1}\beta}{\partial t} + \nabla^{m-1}\beta * \nabla Rm + \nabla^{m-1}\beta * \nabla \gamma\\
& = & \nonumber \nabla \Big( \Delta \nabla^{m-1}\beta + \sum_{i+j = m-1} \nabla^i \beta * \nabla^ j Rm + \sum_{i+j=m} \nabla^{i+1} f * \nabla^j \beta + \sum_{i+j = m-2}\nabla^i \beta * \nabla^{j+1}\gamma\\
&& \nonumber + \sum_{i+j+k = m} \nabla^i (e^{\frac{p+1}{2}f}) * \nabla^j\Psi * \nabla^k \Psi  \Big) + \nabla^{m-1}\beta * \nabla Rm + \nabla^{m-1}\beta * \nabla \gamma
\eea
from which the equation (\ref{eqn:grad beta m}) follows by expanding the terms in the bracket. Note that we need the formula $\nabla \Delta \nabla^{m-1}\beta = \Delta \nabla^m \beta + \nabla^{m-1}\beta * \nabla Rm + \nabla^m\beta * Rm$, which follows from the Riccati equations.
\hfill\qed
\begin{lemma}
We have the evolution equation for $\nabla^m \Psi$
\bea\label{eqn:grad psi m}
\frac{\partial \nabla^m\Psi}{\partial t}& = & \Delta  \nabla^m \Psi + \sum_{i+j = m} \nabla^i \Psi * \nabla^j Rm + \sum_{i+j=m+1} \nabla^{i+1} f * \nabla^j \beta \\
&& \nonumber + \sum_{i+j+k = m+1} \nabla^i(e^{\frac{p+1}{2}f}) * \nabla^j \Psi * \nabla^k \Psi + \sum_{i+j = m-1} \nabla^i \Psi * \nabla^{j+1}\gamma.
\eea
\end{lemma}
\noindent{\it Proof. } When $m=1$, this equation is given by (\ref{eqn:grad psi 1}). The general equation follows by similar calculation as  in deriving the equation (\ref{eqn:grad beta m}) for $\nabla^m \beta$, so we omit the details.\hfill\qed

\begin{lemma}
The evolution equation for $\nabla^{m+1} f$ is given by
\bea\label{eqn:grad f m}
\frac{\partial \nabla^{m+1}f}{\partial t}& = & \Delta \nabla^{m+1} f + \sum_{i+j = m-1} \nabla^{i+2} f * \nabla^j Rm + \nabla^{m+1} e^{-f}\\
&&\nonumber + \sum_{i+j+k=m+1} \nabla^i (e^{-(p+1)f}) * \nabla ^j \beta * \nabla^k \beta + \sum_{i+j=m+1} \nabla^i\Psi * \nabla^j \Psi\\
&& \nonumber + \sum_{i+j = m-1} \nabla^{i+1}f * \nabla^{j+1}\gamma.
\eea
\end{lemma}
\noindent{\it Proof. } When $m=1$, this equation is given by (\ref{eqn:nabla2f 1}). The general equation can be proved similarly as in Lemma \ref{lemma beta}.\hfill\qed

\begin{lemma}
The evolution equation for $\nabla^m Rm$ is
\bea\label{eqn:grad Rm m}
\frac{\partial\nabla^{m-1} Rm}{\partial t}& = & \Delta \nabla^{m-1} Rm + \sum_{i+j = m-1} \nabla^i Rm * \nabla^j Rm + \sum_{i+j=m-1} \nabla^{i+2 } f * \nabla^{j+2} f\\
&&\nonumber + \sum_{i+j+k+l = m-1}\nabla^i Rm* \nabla^j (e^{-(p+1)f}) * \nabla^k \beta * \nabla^l \beta +\sum_{i+j=m+1}\nabla^i \Psi*\nabla^j\Psi \\
&&\nonumber   + \sum_{i+j+k = m-1} \nabla^i Rm * \nabla^j \Psi * \nabla^k \Psi + \sum_{i+j+k=m+1} \nabla^i (e^{-(p+1)f} ) * \nabla^j\beta * \nabla^k \beta\\
&&\nonumber + \sum_{i+j = m-2} \nabla^i Rm * \nabla^{j+1}\gamma.
\eea
\end{lemma}
\noindent{\it Proof. } When $m=2$ the equation is (\ref{eqn:grad Rm 1}). The general case follows by induction similar as in Lemma \ref{lemma beta}. \hfill \qed

\medskip

\noindent{\it Proof of Theorem \ref{prop:4}. } We will use induction to prove (\ref{eqn:grad m}). When $m=2$, this is proved in Lemma \ref{prop:3}. So we assume the estimate has been proved for any nonnegative integer  no bigger than $ m-1$. Our goal is to prove (\ref{eqn:grad m}) for $m$. By induction assumption, the following hold for a constant $C=C(m,K,T)>0$: $\sup_{M^{10-p}} |\nabla f|\le \frac{C}{t^{1/2}}$, and
\begin{equation}\label{eqn:my 1}\sup_{M^{10-p}}( |\nabla^k\beta| + |\nabla^k\Psi| + |\nabla^{k-1} Rm| + |\nabla^{k+1} f| ) \le \frac{C}{t^{k/2}},\quad \forall ~ k\le m-1,\end{equation}
%$$\sup_{M^{10-p}} |\nabla^k Rm|\le \frac{C}{t^{\frac{k+1}{2}}},\quad \forall ~ k\le m-2,$$
%$$\sup_{M^{10-p}} |\nabla^k f| \le \frac{C}{t^}$$
\begin{equation}\label{eqn:my 2}\sup_{M^{10-p}}\Big( |\nabla^i e^{\frac{p+1}{2}f}| + |\nabla^i e^{- (p+1)f}| + |\nabla^i e^{-f}|\Big)  \le C\Big(\frac{1}{t^{i/2}} + 1  \Big),\quad \forall ~ i\le m. \end{equation}
\begin{equation}\label{eqn:my 3}\sup_{M^{10-p}}|\nabla^{j+1} \gamma|\le C \Big( 1+ \frac{1}{t^{(j+1)/2}}  \Big),\quad \forall ~ j\le m-1.\end{equation}
We denote $$G_m: = |\nabla^{m-1} Rm|^2 + |\nabla^m \beta|^2 + |\nabla^m \Psi|^2 + |\nabla^{m+1} f|^2.$$ We calculate making use of (\ref{eqn:grad beta m}), (\ref{eqn:grad psi m}), (\ref{eqn:grad f m}) and (\ref{eqn:grad Rm m})
\bea\nonumber
&& \ho G_m\\
& \le &\nonumber -2 |\nabla^{m+1} \beta|^2 - 2 |\nabla^{m+1}\Psi|^2 - 2 |\nabla^{m+2} f|^2 - 2 |\nabla^m Rm|^2  \\
&& \nonumber + C|\nabla^m \beta| \Big( \sum_{i+j = m} |\nabla^i \beta || \nabla^ j Rm| + \sum_{i+j=m+1} |\nabla^{i+1} f| | \nabla^j \beta| + \sum_{i+j = m-1}|\nabla^i \beta| | \nabla^{j+1}\gamma|\\
&& \nonumber + \sum_{i+j+k = m+1} | \nabla^i (e^{\frac{p+1}{2}f})|| \nabla^j\Psi| | \nabla^k \Psi| \Big) + C |\nabla ^m \Psi|  \Big( \sum_{i+j = m} | \nabla^i \Psi | | \nabla^j Rm |\\
&&\nonumber  + \sum_{i+j=m+1} |\nabla^{i+1} f || \nabla^j \beta| + \sum_{i+j+k = m+1} | \nabla^i(e^{\frac{p+1}{2}f}) || \nabla^j \Psi|| \nabla^k \Psi| + \sum_{i+j = m-1}  |\nabla^i \Psi| | \nabla^{j+1}\gamma|\Big)\\
%&& \nonumber \\
&&\nonumber+ C  |\nabla^{m+1}f|\Big(  \sum_{i+j = m-1} |\nabla^{i+2} f | | \nabla^j Rm | + |\nabla^{m+1} e^{-f}| + \sum_{i+j=m+1} | \nabla^i\Psi | | \nabla^j \Psi|\\
&&\nonumber + \sum_{i+j+k=m+1} | \nabla^i (e^{-(p+1)f}) | | \nabla ^j \beta | | \nabla^k \beta| + \sum_{i+j = m-1}  | \nabla^{i+1}f | | \nabla^{j+1}\gamma|\Big)\\
%&& \nonumber \\
&& \nonumber+ C |\nabla^{m-1} Rm|  \Big( \sum_{i+j = m-1} | \nabla^i Rm || \nabla^j Rm|  + \sum_{i+j=m-1} |\nabla^{i+2 } f | |  \nabla^{j+2} f|  +\sum_{i+j=m+1}  |\nabla^i \Psi|  | \nabla^j\Psi |\\
&&\nonumber + \sum_{i+j+k+l = m-1} | \nabla^i Rm|  | \nabla^j (e^{-(p+1)f}) | | \nabla^k \beta | | \nabla^l \beta| + \sum_{i+j+k = m-1} | \nabla^i Rm | |  \nabla^j \Psi | | \nabla^k \Psi| \\
&&\nonumber   + \sum_{i+j+k=m+1} | \nabla^i (e^{-(p+1)f} ) | | \nabla^j\beta | | \nabla^k \beta| + \sum_{i+j = m-2} |\nabla^i Rm | |  \nabla^{j+1}\gamma|\Big)
%&&\nonumber
\eea
%\newpage
We will  use (\ref{eqn:my 1}), (\ref{eqn:my 2}) and (\ref{eqn:my 3}) to estimate the terms on the RHS of the above inequality.

\medskip

\noindent $\bullet$ The terms $C |\nabla ^m \beta| (\cdots )$:
\bea\nonumber
&& C|\nabla^m \beta| (\cdots)\\ \nonumber & \le & C|\nabla^m \beta|\Big\{t^{-1/2} ( |\nabla^{m-1} Rm| + |\nabla^{m+1} \beta| + |\nabla^m \beta| + |\nabla^m \Psi| ) +   |\nabla^m \beta|\\
&&\nonumber + |\nabla^m Rm| +  |\nabla^{m+2} f| + |\nabla^{m+1} f| + |\nabla^{m+1}\Psi| + \frac{1}{t^{(m+1)/2}} + \frac{1}{t^{m/2}}   \Big\}\\
& \le & \nonumber  C t^{-1/2} ( |\nabla^m \beta|^2 + |\nabla^{m-1}Rm|^2 + |\nabla^m \Psi|^2  ) + C(1+ t^{-1}) |\nabla^m\beta |^2 + C |\nabla^{m+1} f|^2  \\
&&\nonumber + \frac{1}{10} \Big( |\nabla^{m+1}\beta|^2 + |\nabla ^m Rm|^2 + |\nabla^{m+2} f|^2 + |\nabla^{m+1}\Psi|^2   \Big) + C t^{-(m+1)/2 } |\nabla^m \beta| + C
\eea

\noindent$\bullet$ The terms $C|\nabla^m \Psi|(\cdots)$:
\bea\nonumber
&& C|\nabla^m \Psi|(\cdots)\\ \nonumber & \le & C|\nabla^m \Psi| \Big\{t^{-1/2} ( |\nabla^{m-1} Rm| + |\nabla^{m+1} f| + |\nabla^{m+1} \beta| + |\nabla^m \beta| + |\nabla^m \Psi|   ) + |\nabla^m Rm| \\
&&\nonumber + |\nabla^m \Psi| + |\nabla^{m+2} f| + |\nabla^{m+1} f| + |\nabla ^{m+1}\Psi| + \frac{1}{t^{m/2}} + \frac{1}{t^{(m+1)/2}} + 1 \Big\}\\
&\le & C t^{-1/2} ( |\nabla^m \Psi|^2 + |\nabla^{m-1} Rm|^2 + |\nabla^{m+1} f|^2 +|\nabla^m \beta|^2   ) + C(1+ t^{-1}) |\nabla^m \Psi|^2 \nonumber \\
&& \nonumber +C  |\nabla^{m+1} f|^2  + \frac{1}{10} \Big( |\nabla^{m+1}\beta|^2 + |\nabla^{m+2} f|^2 + |\nabla^m Rm|^2 + |\nabla^{m+1}\Psi|^2   \Big) + \frac{C |\nabla^m \Psi|}{t^{(m+1)/2}} + C.
\eea
\noindent$\bullet$ The terms $C|\nabla^{m+1} f|(\cdots)$
\bea\nonumber
&& C|\nabla^{m+1} f| (\cdots)\\
&\le  & \nonumber C |\nabla^{m+1} f| \Big\{  t^{-1/2} ( |\nabla^{m-1} Rm| + |\nabla^m \Psi| + |\nabla^m \beta|  ) + |\nabla^{m+1} f| + |\nabla^{m+1}\Psi| \\
&& \nonumber + |\nabla^{m+1 }\beta| + \frac{1}{t^{m/2}} + \frac{1}{t^{(m+1)/2}}   \Big\}\\
& \le & \nonumber C t^{-1/2} ( |\nabla^{m+1}f|^2 + |\nabla^{m-1}Rm|^2 + |\nabla^m \Psi|^2 + |\nabla^m \beta|^2  ) + C|\nabla^{m+1} f|^2 \\
&&\nonumber + \frac{1}{10}\Big( |\nabla^{m+1}\Psi|^2 + |\nabla^{m+1}\beta|^2 \Big) + C \frac{|\nabla^{m+1} f|}{t^{(m+1)/2}} + C.
\eea

\noindent$\bullet$ The terms $C |\nabla^{m-1} Rm| (\cdots)$:
\bea\nonumber
&& C |\nabla^{m-1}Rm| (\cdots)\\
& \le & \nonumber C |\nabla^{m-1} Rm| \Big\{ t^{-1/2} ( |\nabla^m \Psi| + |\nabla^m\beta| ) + |\nabla^{m-1} Rm| + |\nabla^{m+1} f| \\
&& \nonumber + |\nabla^{m+1}\Psi | + |\nabla^{m+1}\beta| + \frac{1}{t^{m/2}} + \frac{1}{t^{(m+1)/2}}   \Big\}\\
& \le &\nonumber C t^{-1/2} ( |\nabla^{m-1} Rm|^2 + |\nabla^m\Psi|^2 + |\nabla^m \beta|^2   ) + C |\nabla^{m-1} Rm|^2 + C|\nabla^{m+1} f|^2 \\
&& \nonumber + \frac{1}{10}\Big(   |\nabla^{m+1} \beta|^2 + |\nabla^{m+1}\Psi|^2   \Big) + C \frac{|\nabla^{m-1} Rm|}{t^{(m+1)/2}} + C.
\eea
Combining the inequalities above, we get
\bea\nonumber
&& \ho ( t^m G_m) \\
& = &\nonumber m t^{m-1} G_m + t^m \ho G_m \\
& \le &\nonumber m t^{m-1} G_m + t^m \Big\{ -|\nabla^{m+1} \beta|^2 - |\nabla^{m+1}\Psi|^2 - |\nabla^{m+2} f|^2 - |\nabla^m Rm|^2 \\
&& \nonumber + C t^{-1/2} G_m + C G_m + C t^{-1}  ( |\nabla^m\beta|^2 + |\nabla^m \Psi|^2) + C\Big\} + C t^{(m-1)/2} \Big(|\nabla^m\beta|\\
&&\nonumber + |\nabla^m \Psi| + |\nabla^{m+1} f| + |\nabla^{m-1}Rm|  \Big)\\
&\le & C t^{m-1} G_m  -  t^m G_{m+1} + C,\label{eqn:induction}
\eea
where in the last step we applied Cauchy-Schwarz inequality. Note that from the calculations above it is not hard to see that the inequality (\ref{eqn:induction}) in fact holds for $t^i G_i$ for any $1\le i\le m$, maybe with different constants which depend only on $K, T, m$. Define
$$H: = t^m G_m + \sum_{i=1}^{m-1} B_i t^i G_i + B_0 ( |\nabla f|^2 + |\beta|^2 + |\Psi|^2 ).$$
For suitable choice of the constants $B_i$'s, we have
\bea\nonumber
\ho H \le - t^m G_{m+1} + C,
\eea for some constant $C = C(K,T,m)>0$. From maximum principle, it follows that $\sup_{M^{10-p} } H \le C(K,T,m)$, thus $\sup_{M^{10-p}} G_m \le \frac{C(K,T,m)}{t^m}$ for any $t\in (0,T)$.  Thus we finish the proof of Theorem \ref{prop:4}.

\bigskip
Once we have the estimates of Theorem \ref{prop:4}, the proof for part (c) of Theorem \ref{eqn:system} can be completed in the same way as for part (c) of Theorem \ref{Euclidean}. \hfill \qed

\bigskip
We conclude by observing that Theorem \ref{eqn:system} implies similar theorems for the simpler dimensional reductions described in sections \ref{section2.2} and \ref{section2.3}.
In some cases, under suitable assumptions, the theorem can be strengthened. For example, in the case \S 2.3 and if $\tilde\lambda\leq 0$, since $\beta$ is automatically $0$, it is easy to obtain a bound on $|f|$ on any finite time interval $[0,T)$. Thus the maximum time $T$ of existence is now characterized by
\bea
\limsup_{t\to T^-} \sup_{M^{10-p}} ( |Rm|  +  |\Psi|   ) = \infty.
\eea
Also, when the dimension $p$ satisfies $p\geq 4$, the forms $\beta$ and $\Psi\wedge\Psi$ are automatically $0$, and the equations simplify a great deal.

\begin{appendix}

\section{Conventions}
\setcounter{equation}{0}

We consider metrics $g$ on an $11$-dimensional manifold $M^{11}$, which can be either Lorentzian or Riemannian. We distinguish between the two cases by a number $\sigma$, which is defined to be $+1$ if $M^{11}$ is Lorentzian, and $-1$ if $M^{11}$ is Riemannian.

\medskip
The $4$-form $F$ is expressed in components as
\bea\nonumber
F={1\over 4!}F_{ABCD}\,dx^A\wedge dx^B\wedge dx^C\wedge dx^D
\eea
and $F^2_{AB}$, $|F|^2$ and $|\nabla F|^2$ are defined by
\bea\nonumber
F^2_{AB}={1\over 3!}F_{ACDE}F_B{}^{CDE},
\quad
|F|^2={1\over 4!}F_{ACDE}F^{ACDE},
\quad
|\nabla F|^2={1\over 4!} F_{ACDE,B}F^{ACDE,B}.
\eea

The operator $\Box_g = d d^\dagger + d^\dagger d$ is the Hodge-Laplacian. We also need the Laplacian $\Delta_g$ defined on tensors or forms by
$\Delta_g=g^{AB}\nabla_B\nabla_A$. The Lichnerowicz-Weitzenb\"ock formula says that the two differ by curvature terms, e.g. on forms such as $F$,
\bea\nonumber
\Box_gF=-\Delta_gF+Rm * F
\eea
where $Rm * F$ denotes a pointwise bilinear expression in the components of $Rm $ and $F$.

\smallskip

The Hodge $\star$ is defined by
\bea\nonumber
\alpha\wedge \star\beta=\<\alpha,\beta\>\sqrt{-\sigma g}\,d^{11}x
\eea
To emphasize the metric sometimes we write $\star_g$ instead of $\star$.
In particular, if we set $d vol_g=\sqrt{-\sigma g}d^{11}x$, then $\star dvol_g=-\sigmaone $ and $\star 1=dvol_g$. The adjoint $d^\dagger$ of the exterior derivative $d$ acting on 4-forms is then given by
\bea
d^\dagger=-\sigmaone \star d\star.
\eea

\section{The field equations of $11D$ supergravity}
\setcounter{equation}{0}

We provide here a derivation of the field equations of $11D$ supergravity for the convenience of the reader. Let $g$ be a metric on $M^{11}$ and $A$ be a $3$-form. Let $F= dA$. Recall that we set $\sigma=+1$ if $g$ is Lorentzian, and $\sigma=-1$ if $g$ is Riemannian. The action of $11D$ supergravity is given by
\begin{equation}\label{eqn:functional}
{\mathcal L}[g, A] = \int_M R \sqrt{-\sigma g} - {1\over 2} F \wedge \star F + {1\over 6} F \wedge F \wedge A,
\end{equation}
The Euler-Langrange equations of the functional $\mathcal L$ in (\ref{eqn:functional}) are readily found to be
\bea
\label{eqn:app1}
d\star F  &= & \frac 1 2 F \wedge F \nonumber \\
R_{AB} - \frac 1 2 R g_{AB} &=& {1\over 2}F^2_{AB} - \frac 1 4 |F|^2 g_{AB}.
\eea
Contracting the second equation in (\ref{eqn:app1}), we see that $R = \frac 1 6 |F|^2$. Therefore we find that this equation is equivalent to
\begin{equation}\label{eqn:app2}
R_{AB} = \frac 1 2 F^2_{AB}- \frac 1 6 |F|^2 g_{AB}.
\end{equation}
This is the form of the equation used in (\ref{FieldEqs}).
Moreover, if we apply $d\star$ on both sides of the first equation in (\ref{eqn:app1}), we get
$$dd^\dagger F = - \frac \sigma 2 d\star (F\wedge F),$$
combining with $dF = 0$ we see $\Box_g F = -  \frac \sigma 2 d\star(F\wedge F)$.

\section{Field equations and stationary points}
\setcounter{equation}{0}

Clearly the solutions of the field equations must be stationary points of the dynamical system. We discuss here some simple situations when the converse is true.

\medskip
Suppose that the cohomology group $H^3(M^{11},{\bf R})$ is $0$. Then at a stationary point, we have $\star d\star F-{1\over 2}\star (F\wedge F)=dB$ for some smooth $2$-form $B$. Applying $d\star$ to both sides of this equation and using the fact that $F$ is closed gives $\star d\star dB=0$, i.e. $d^\dagger dB=0$.

If $M^{11}$ is compact, this implies that $|dB|^2=0$, and if $M^{11}$ is also Riemannian, this implies that $dB=0$. The field equations are then satisfied.
Thus the field equations (\ref{FieldEqs}) are equivalent to the stationary point condition of the flow (\ref{11Dflow}) when $M^{11}$ is a compact Riemannian manifold with vanishing cohomology group $H^3(M^{11},{\bf R})$.

\medskip
More generally, it follows from the stationary point condition that the form $\alpha$
\bea
\alpha=\star d\star F-{1\over 2}\star (F\wedge F)
\eea
is both closed and co-closed, i.e., $d\alpha=0$ and $d^\dagger\alpha=0$. Thus $\alpha$ must be $0$ and the field equations are satisfied if $M^{11}$ is assumed not to have any such form which is non-trivial. This is equivalent to the non-existence of non-trivial harmonic forms when $M^{11}$ is Riemannian.

\end{appendix}

%\newpage

\bigskip
Department of Mathematics, Columbia University, New York, NY 10027

tfei@math.columbia.edu, bguo@math.columbia.edu,
phong@math.columbia.edu

\end{document}